\documentclass[12pt, oneside,reqno]{birkjour}
\usepackage{amsfonts}
\usepackage{amsmath, amssymb}
\usepackage{graphicx}
\usepackage{subfigure}
\usepackage[font=small,labelfont=bf]{caption}
\usepackage{epstopdf}
\usepackage[pdfpagelabels,hyperindex]{hyperref}
\usepackage{xcolor}
\usepackage{amsthm}
\usepackage{comment}
\usepackage{bbm}
\usepackage{float}
\usepackage{pgfplots}
\usepackage{listings}
\usepackage{longtable}
\usepackage{mathrsfs}
\usepackage{enumerate}
\hypersetup{
pdftitle={},
pdfsubject={},
pdfauthor={},
pdfkeywords={}
}

 \newtheorem{thm}{Theorem}[section]

 \newtheorem{prop}[thm]{Proposition}
 \theoremstyle{definition}
 \newtheorem{definition}[thm]{Definition}
 \theoremstyle{remark}
 \newtheorem{rmk}[thm]{Remark}
 \newtheorem{exmp}{Example}

\begin{document}

%
%
%
%
%
%

\title[]{ A Semidiscrete Lagrangian-Eulerian scheme for the LWR traffic model with discontinuous flux}

\author[E. Abreu]{Eduardo Abreu}
\address{%
Universidade Estadual de Campinas (UNICAMP), 
Department of Applied Mathematics, IMECC, 13083-859, 
Brazil.}
\email{eabreu@unicamp.br}

\author[M. T. Chiri]{Maria Teresa Chiri}
\address{Department of Mathematics and Statistics, Queen’s University, Kingston, ON K7L3N6, Canada}
\email{maria.chiri@queensu.ca}

\author[R. De la cruz]{Richard De la cruz}
\address{Universidad Pedagógica y Tecnológica de 
Colombia (UPTC), School of Mathematics and Statistics, 150003, 
Colombia.}
\email{richard.delacruz@uptc.edu.co}

\author[J. Juajibioy]{Juan  Juajibioy}
\address{Universidad Pedagógica y Tecnológica de 
Colombia (UPTC), School of Mathematics and Statistics, 150003, 
Colombia.}
\email{juan.juajibioy@uptc.edu.co}

\author[W. Lambert]{Wanderson Lambert}
\address{ICT - UNIFAL, Cidade Universitária - 
BR 267 Km 533 Rodovia José Aurélio Vilela, Poços de Caldas, 
Alfenas, MG, Brazil.}
\email{wanderson.lambert@unifal-mg.edu.br}

\subjclass{Primary 35L50, 35L65,35R05; Secondary 76S05}

\keywords{Conservation laws with discontinuous flux, semidiscrete Lagrangian-Eulerian scheme}

\date{today}
\begin{abstract}
 {In this work, we present a semi-discrete scheme to approximate solutions to the scalar LWR traffic model with spatially discontinuous flux, described by the equation  $u_t + (k(x)u(1-u))_x = 0$.} This approach is based on the Lagrangian-Eulerian method proposed by E. Abreu, J. Francois, W. Lambert, and J. Perez \emph{[J. Comp. Appl. Math. 406 (2022) 114011]} for scalar conservation laws. 
 We derive a non-uniform bound on the growth rate of the total variation for approximate solutions. Since the total variation can explode only at $x=0$, we can provide a convergence proof for our scheme in $BV_{loc}(\mathbb{R}\setminus \lbrace 0 \rbrace)$ by using Helly's compactness theorem.

\end{abstract}
\maketitle

\section{Introduction}\label{Intro}
In this work, we introduce and analyze a semi-discrete scheme, 
in the Lagrangian-Eulerian framework, to compute approximate 
solutions to the Cauchy problem for the following class of scalar 
conservation laws with spatially discontinuous flux:
\begin{align}\label{1-1}
&\partial_t u+\partial_x F(k(x),u)=0, \ (t,x)\in (0,T)\times \mathbb{R}\\
&u(0,x)=u_0(x), \quad x\in \mathbb{R}.\label{1-1-D}
\end{align}
Following \cite{Towers1, Vovelle}, we consider a flux function with a unique 
spatial discontinuity at \( x = 0 \) of the form:
\begin{equation}\label{1-2}
F(x,u)=k(x)f(u),
\end{equation}
where $f:[0,1]\to \mathbb{R}$ is given by
\begin{equation*}
f(u)=u(1-u),
\end{equation*}
and
\begin{equation}\label{k-f}
k(x)=k_L(1-H(x))+ k_RH(x).
\end{equation}
Here, $H$ represents the Heaviside function
\begin{equation*}
H(x)=\begin{cases}
0, \ \text{if} \ x<0,\\
1, \ \text{if} \ x>0,\\
\end{cases}
\end{equation*}
and  $k_L$, $k_R$, is a pair of positive  { constants} such that 
$k_L\neq k_R$. \\
{
Equation \eqref{1-1} can be seen as an LWR traffic model \cite{LW,R} in which a change 
of the flow-density relation occurs at \( x = 0 \). The classic LWR model expresses 
the conservation of the total number of vehicles and postulates that the average traffic 
speed \( v(x,t) \) is a function of the traffic density \( u(t,x) \) alone. Thus, the mean 
traffic flow (the number of cars crossing the point \( x \) per unit of time) is given by  $
f(x,t) = u(x, t) v(u(x, t)),$
and we are led to the hyperbolic conservation law 
$
u_t + (uv(u))_x = 0.
$
In the case of Equation \eqref{1-1}, the speed function is given by:
$$ v(u(x, t))=\begin{cases} k_L(1-u) \quad &\text{if } x<0\\
 k_R(1-u)\quad &\text{if } x>0\end{cases}, $$
hence, we are modeling traffic on two stretches of road of infinite length and the same capacity, but with different maximal speeds.
\\

More generally, starting from the seminal works of Isacson \& Temple \cite{Tem} and Risebro et al.\ \cite{Gimse1, Klausen}, conservation laws with discontinuous flux have been extensively studied over the last three decades. The significant interest in this type of partial differential equation is motivated by their application in many different areas of physics and engineering, including two-phase flow models in porous media with changing rock types (for oil reservoir simulation) \cite{Gimse1, Gimse2}; slow erosion granular flow models \cite{Shen2}; Saint Venant models of blood flow in endovascular treatments \cite{Canic, Quarteroni}; and traffic flow models with roads of varying amplitudes or surface conditions \cite{Mochon}.\\

From the mathematical standpoint, the presence of a unique spatial discontinuity at \( x = 0 \) in the flux \eqref{1-2} is already sufficient to significantly affect the behavior of solutions of \eqref{1-1}, indeed:
\begin{enumerate}[(1)]
\item \textit{\textbf{Interface entropy conditions are needed to select a unique solution.}}\\
Equation \eqref{1-1} is usually supplemented with appropriate coupling conditions imposed at the point of discontinuity of the flux to guarantee the uniqueness of solutions to the Cauchy problem. Moreover, various types of admissibility conditions (interface entropy conditions) imposed on the left and right traces of the solution have been introduced in the literature, according to different modeling assumptions \cite{Andreianov1, Andreianov2}. These conditions lead to different solutions of the Cauchy problem \eqref{1-1}, \eqref{1-1-D}.
\item \textit{\textbf{The discontinuity in the flux can increase the total variation.}}\\
This may happen regardless of the choice of interface entropy condition. To illustrate this phenomenon, we follow \cite{mishra1} and consider the case of \( k(x) = 1 \cdot (1-H(x)) + 2 \cdot H(x) \) with the initial datum \( u_0(x) = 0.5 \), whose total variation is clearly equal to zero.  
By using the Riemann solver \cite{Mochon} introduced by Mochon in this setting (which coincides with the AB-entropy condition of \cite{Adimurthi1} with \( A = \frac{1}{2} \) and \( B = \frac{1}{2} - \frac{1}{2\sqrt{2}} \)), we find the solution:
\begin{equation}\label{misrha-ex1}
u(t,x)=\begin{cases}
\frac{1}{2}, \ \text{if}\ x<0,\\
\frac{1}{2}-\frac{1}{2\sqrt{2}}, \ \text{if}\ 0<x<\sigma t,\\
\frac{1}{2}, \ \text{if}\ x>\sigma t,
\end{cases}
\end{equation}
where $\sigma$ is the classical Rankine-Hugoniot relation applied 
to the flux $F=2u(1-u)$ with values $u_-=-\frac{1}{2\sqrt{2}}+\frac{1}{2}$, $u_+=\frac{1}{2}$.
 Notice that in this case
\[
TV(u(t,\cdot))> 0.
\]
\end{enumerate} 
The general lack of uniform bounds on the total variation of entropy solutions when the flux presents spatial discontinuity constitutes an obstruction to the application of the TVD (Total Variation Diminishing) methodology. Therefore, several efforts have been made to provide BV regularity of entropy solutions by adding constraints to the shape of the flux. For instance, in \cite{Piccoli}, the authors consider a class of flux functions for which there exists a pair of positive constants \( c \) and \( C \) such that
\begin{equation}\label{Lipflux}
c|u-\overline{u}|\leq |F(x,u)-F(x,\overline{u})|\leq C|u-\overline{u}|.
\end{equation}
This property has also been used in \cite{Towers3} to provide bounds on the total variation of the numerical flux functions constructed by finite difference schemes and to obtain a uniform estimate of the total variation of approximate solutions. However, both \cite{Piccoli} and \cite{Towers3} also require strict monotonicity of the flux, which is not a suitable assumption in traffic modeling and indeed does not apply in our case. In this context, we note that several efforts have been made to construct second-order schemes that satisfy the flux total variation diminishing property (FTVD) \cite{Kenneth4, Adimurthi}.
 Finally, for semi-discrete approximations based on a central upwind scheme, we refer to \cite{Wang}.
\\

We recall that, in the context of scalar conservation laws with discontinuous flux, the notion of a weak solution agrees with the classical one, i.e., following \cite{Bressan}
\begin{definition}[\textbf{Weak solutions}]
A function \( u \in L^1_{\text{loc}}(\mathbb{R}_+ \times \mathbb{R}) \) 
is a \textit{weak solution} to the Cauchy problem \eqref{1-1}-\eqref{1-1-D} if it satisfies
\begin{equation*}
\int_{0}^{\infty}\int_{\mathbb{R}}\Big(u(t,x)\partial_t \phi(t,x)+F(k(x),u(t,x))\partial_x\phi(t,x)\Big)dtdx+\int_{\mathbb{R}}u_0(x)\phi(0,x)dx=0,
\end{equation*}
for all $\phi\in C^{\infty}_{c}(\mathbb{R}_{+}\times \mathbb{R})$.
\end{definition}
If \( u \) is a weak solution of the Cauchy problem \eqref{1-1}, \eqref{1-1-D},  
then, from \cite[Section 3]{mishra1}, \( u \) satisfies the following 
Rankine-Hugoniot relations:
\begin{equation}\label{4-lok-RH}
f(u^{-}(t))=f(u^+(t))
\end{equation}
where $u^{-}(t)=\lim_{x\to 0^-}u(t,x)$ and $u^{+}(t)=\lim_{x\to 0^+}u(t,x)$.\\
To select a unique solution, we consider the notion of generalized entropy used in \cite{Vovelle, Towers1}, i.e., for the entropy flux:
\begin{equation*}
\Phi(u(t,x),\kappa)=\text{sgn}\big(u(t,x)-\kappa\big)(f(u(t,x))-f(\kappa)).
\end{equation*}
we say that
\begin{definition}[\textbf{Entropy solution}]
A bounded and measurable function \( u : \mathbb{R}_+ \times \mathbb{R} \to \mathbb{R} \) 
is a \textit{weak entropy solution} of the Cauchy problem \eqref{1-1}, \eqref{1-1-D} if
\begin{align}
&\int_{0}^{\infty}\int_{\mathbb{R}}\Big(|u(t,x)-\kappa|\partial_t \phi(t,x)+k(x)\Phi(u(t,x),\kappa)\partial_x\phi(t,x)\Big)dt\,dx \notag\\
&+\int_{\mathbb{R}}{|u_0(x)-\kappa|}\phi(0,x)dx+|k_{R}-k_{L}|f(\kappa)\int_{0}^{\infty}\phi(t,0)dt\geq 0, \label{3-loc-weak}
\end{align}
for all $\kappa \in [0,1]$ and for all $\phi\in C^{\infty}_{c}(\mathbb{R}_+\times \mathbb{R})$.
\end{definition}

{ This notion of entropy was originally derived by considering the problem \eqref{1-1}-\eqref{1-1-D} with a continuous approximation of the flux. We point out that it effectively selects a unique solution because of the strict concavity of the flux outside the discontinuity. Without this hypothesis, a further crossing condition of the right and left fluxes would be needed to ensure uniqueness \cite{Andreianov}. Moreover, the solution selected by this entropy condition coincides with the one determined by the critical \( AB \)-connection \cite{Adimurthi1}.
 \\

 Several papers have introduced different Lagrangian-Eulerian methods (see, e.g., \cite{Doug1} for convection-diffusion problems, \cite{Aq} for linear transport problems, and \cite{Doug2} for nonlinear scalar transport problems). These formulations are the result of efforts to develop fast, accurate, and stable versions of the modified method of characteristics \cite{Doug2, Huang}. The semidiscrete Lagrangian-Eulerian scheme developed in this work approximates the entropy solution for \eqref{1-1} selected by \eqref{3-loc-weak}. As far as we know, this is the first work where Lagrangian-Eulerian schemes are employed within the context of conservation laws with discontinuous flux.
}
\\

This paper is organized as follows: In Section \ref{sdlea}, we construct the semidiscrete 
Lagrangian-Eulerian scheme for the LWR traffic model with discontinuous 
flux. We perform the convergence analysis of the scheme, including the corresponding \( L^{\infty} \)-estimates and \( L^{1} \)-estimates, as well as the non-uniform TV-estimates and temporal continuity, in Section \ref{aansd}. The convergence to the weak entropy unique solution is addressed in Section \ref{cwesol}. The numerical experiments are presented and discussed in Section \ref{numexa}. Finally, our concluding remarks and future research directions are presented in Section \ref{ConRem}.

\section{A semidiscrete Lagrangian-Eulerian approach}\label{sdlea}
In this section, we design and analyze a new semi-discrete scheme to construct approximate solutions for the scalar conservation law with discontinuous flux introduced in \eqref{1-1}. This approximation is based on the concept of the no-flow property introduced in \cite{Abreu-John} and applied to design semidiscrete schemes for scalar conservation laws in \cite{Eduardo1}.
\subsection{A semi-discrete formulation}\label{Subsec-2.1}
Given a time horizon \( T > 0 \), we discretize the equation \eqref{1-1} 
over the time-space region \([0, T] \times \mathbb{R}\) by selecting 
uniform mesh size parameters \( \Delta x > 0 \) and \( \Delta t > 0 \). 
Let \( x_{j}^n = j \Delta x \) and 
\[
x^n_{j \pm \frac{1}{2}} = \left( j \pm \frac{1}{2} \right)\Delta x,
\]
for \( j = 0, \pm 1, \pm 2, \dots \), \( t^n = n \Delta t \), \( n \in \mathbb{N} \), we define
 
\begin{align} 
& u_{j+\frac{1}{2}}^n=\frac{1}{\Delta x}\int_{x^n_{j}}^{x^{n}_{j+1}}u(t^n,x)dx,\notag\\
&D_j^n=\left\{(x,t):\text{ for } \sigma_{j-\frac{1}{2}}(t)\leq x
\leq \sigma_{j+\frac{1}{2}}(t), \ t^n\leq t\leq t^{n+1}\right\},\label{Djn}
\end{align}
where the curves $\sigma_{j\pm\frac{1}{2}}(t)$ enclosing the region $D_j^n$ are called {\it no-flow curves} and are solutions to the initial value problem (IVP)
\begin{align}
&\frac{d}{dt}\sigma_{j+\frac{1}{2}}(t)=\frac{F\big(k(\sigma_{j+\frac{1}{2}}(t)),u(t,\sigma_{j+\frac{1}{2}}(t))\big)}{u(t,\sigma_{j+\frac{1}{2}}(t))},\label{ODE-sigma}\\
&\sigma_{j+\frac{1}{2}}(t^n)=x^n_{j+\frac{1}{2}}.\notag
\end{align}
{Consider} a Taylor's first order approximation of   $\sigma_{j+\frac{1}{2}}(t)$, 
by replacing the first  derivative with \eqref{ODE-sigma}, we obtain 
the following approximation of  $\sigma_{j+\frac{1}{2}}(t^{n+1})$ that we indicate as $\overline{x}^n_{j+\frac{1}{2}}$:
\begin{equation*}
\overline{x}^n_{j+\frac{1}{2}}=x^{n}_{j+\frac{1}{2}}+\Delta t\frac{F(k(x^n_{j+\frac{1}{2}}),u^n_{j+\frac{1}{2}})}{u^n_{j+\frac{1}{2}}}=\begin{cases}
x^{n}_{j+\frac{1}{2}}+\Delta t k_L(1-u^n_{j+\frac{1}{2}}), \ \text{if} \ j\leq -1,\\
x^{n}_{j+\frac{1}{2}}+\Delta t k_R(1-u^n_{j+\frac{1}{2}}), \ \text{if} \ j\geq 0.\\
\end{cases}
\end{equation*}
In Figure \ref{fig:0}, we present a geometric description of these curves obtained in a specific mesh grid. Future interaction of the no-flow curves can be avoided if we consider the following restriction on the mesh parameter:
\[
x^n_j\leq x^n_{j+\frac{1}{2}}+\Delta t\frac{F(k(x^n_{j+\frac{1}{2}}),u^n_{j+\frac{1}{2}})}{u^n_{j+\frac{1}{2}}}\leq x^n_{j+1},
\]
which implies the condition
\begin{equation}\label{CFL-0}
-\frac{1}{2}\frac{\Delta x}{\Delta t}\leq \frac{F(k(x^n_{j+\frac{1}{2}}),u^n_{j+\frac{1}{2}})}{u^n_{j+\frac{1}{2}}}\leq \frac{1}{2}\frac{\Delta x}{\Delta t}.
\end{equation}
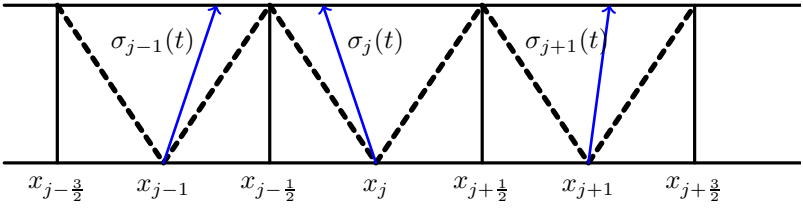
\begin{figure}[h!]
    \centering
\begin{tikzpicture}[scale=0.7,black,line cap=round,line join=round,x=1cm,y=1cm]\clip(-9,-5) rectangle (7,0);
\draw [line width=1.2pt] (-8,-4)-- (7,-4);
\draw [line width=1.2pt] (-7,-4)-- (-7,-1);
\draw [line width=1.2pt] (-3,-4)-- (-3,-1);
\draw [line width=1.2pt] (1,-4)-- (1,-1);
\draw [line width=1.2pt] (5,-4)-- (5,-1);
\draw [line width=2pt, dashed] (-7,-1)-- (-5,-4);
\draw [line width=2pt, dashed] (-5,-4)-- (-3,-1);
\draw [line width=2pt, dashed] (-3,-1)-- (-1,-4);
\draw [line width=2pt, dashed] (-1,-4)-- (1,-1);
\draw [line width=2pt, dashed] (1,-1)-- (3,-4);
\draw [line width=2pt,dashed] (3,-4)-- (5,-1);
\draw [line width=1.2pt] (-8,-1)-- (7,-1);
\draw [->,line width=1pt,blue] (-5,-4) -- (-4,-1);
\draw [->,line width=1pt,blue] (-1,-4) -- (-2,-1);
\draw [->,line width=1pt,blue] (3,-4) -- (3.4,-1);
\draw (-7,-4.5) node[] {$ x_{j-\frac{3}{2}} $};
\draw (-5,-4.5) node[] {$ x_{j-1} $};
\draw (-3,-4.5) node[] {$ x_{j-\frac{1}{2}} $};
\draw (-1,-4.5) node[] {$ x_{j} $};
\draw (1,-4.5) node[] {$ x_{j+\frac{1}{2}} $};
\draw (3,-4.5) node[] {$ x_{j+1} $};
\draw (5,-4.5) node[] {$ x_{j+\frac{3}{2}} $};
\draw (-5.2,-1.7) node[] {$ \sigma_{j-1}(t) $};
\draw (-1,-1.7) node[] {$ \sigma_{j}(t) $};
\draw (2.6,-1.7) node[] {$ \sigma_{j+1}(t) $};
\end{tikzpicture}
    \caption{Behavior of no-flow curves satifaying $|\sigma_{j}(t)|<\frac{1}{2}\frac{\Delta x}{\Delta t}$.}
    \label{fig:0}
\end{figure}
To approximate \( u(t,x) \) at time \( t^{n+1} \), we follow these steps:\\

\textbf{Step 1.} Following \cite{Eduardo1}, at time \( t^n \), we define the linear function:
\begin{align}\label{L-1}
L(x)=\sum_{j}\mathbbm{1}_{I_j}(x)L_{j+\frac{1}{2}}(x),
\end{align}
where
\begin{align*}
L_{j+\frac{1}{2}}(x)=u^n_{j+\frac{1}{2}}+(x-x_{j+\frac{1}{2}})s^n_{j+\frac{1}{2}}, \ x\in I_j=(x_{j},x_{j+1}),
\end{align*}
and $s^n_{j+\frac{1}{2}}$ is defined by
\begin{equation}\label{slope}
s^n_{j+\frac{1}{2}}=\text{minmod}\Big(\frac{u^n_{j+\frac{1}{2}}-u^n_{j-\frac{1}{2}}}{\Delta x},\frac{u^n_{j+\frac{3}{2}}-u^n_{j-\frac{1}{2}}}{2\Delta x},\frac{u^n_{j+\frac{3}{2}}-u^n_{j+\frac{1}{2}}}{\Delta x}\Big)
\end{equation}
with 
\begin{equation}\label{minmod}
\text{minmod}(a_1,a_2,a_3)=
\begin{cases}
s\cdot \text{min}_{i=1,2,3}(|a_i|), \ \text{if}\ s=\text{sgn}(a_1)=\text{sgn}(a_2)=\text{sgn}(a_3).\\
0,\ \text{else}. 
\end{cases}
\end{equation}
Let
\begin{equation*}
\ell^{n+1}_j=\sigma_{j+\frac{1}{2}}(t^{n+1})-\sigma_{j-\frac{1}{2}}(t^{n+1}),
\end{equation*}
we consider the averaged integral of \eqref{1-1} over the region \( D^n_j \). By applying the Divergence theorem and noting that \( \sigma^n_{j}(t) \) is a solution of \eqref{ODE-sigma}, we obtain:
\begin{align}
&\frac{1}{\ell^{n+1}_j}\int_{\sigma_{j-\frac{1}{2}}}^{\sigma_{j+\frac{1}{2}}}u(t^{n+1},x)dx\notag\\
&=\frac{1}{\ell^{n+1}_j}\int_{x_{j-\frac{1}{2}}}^{x_{j+\frac{1}{2}}}u(t^n,x)dx=\frac{1}{\ell^{n+1}_j}\Big(\int_{x_{j-\frac{1}{2}}}^{x_{j}}u(t^n,x)dx+\int_{x_{j}}^{x_{j+\frac{1}{2}}}u(t^n,x)dx\Big).\label{step-1}
\end{align}
At time \( t^n \), at the discontinuity point \( x = 0 \), we apply the Rankine-Hugoniot relation \eqref{4-lok-RH}, \( F(0^-, u(t^n, 0^-)) = F(0^+, u(t^n, 0^+)) \), and obtain:
\begin{align*}
&\frac{1}{\ell^{n+1}_0}\int_{\sigma_{-\frac{1}{2}}}^{\sigma_{\frac{1}{2}}}u(t^{n+1},x)dx+\int_{t^n}^{t^{n+1}}\Big(F(0-,u(t,0-)-F(0+,u(t,0+))\Big)dt\notag\\
&=\frac{1}{\ell^{n+1}_0}\int_{x_{-\frac{1}{2}}}^{x_{\frac{1}{2}}}u(t^n,x)dx=\frac{1}{\ell^{n+1}_0}\Big(\int_{x_{-\frac{1}{2}}}^{0}u(t^n,x)dx+\int_{0}^{x_{\frac{1}{2}}}u(t^n,x)dx\Big).
\end{align*}
Set
\begin{equation}\label{f_j}
f^n_{j+\frac{1}{2}}=\frac{F\big(k(x^n_{j+\frac{1}{2}}),u^n_{j+\frac{1}{2}}\big)}{u^n_{j+\frac{1}{2}}}=k_{j+\frac{1}{2}}(1-u^n_{j+\frac{1}{2}}),
\end{equation}
and
\begin{equation*}
h^n_{j}=\Delta x+\Delta t\big(f^n_{j+\frac{1}{2}}-f^n_{j-\frac{1}{2}}\big),
\end{equation*}
by using the function $L(x)$ defined in \eqref{L-1}, we obtain 
an approximation of \eqref{step-1} as follows
\begin{align}
u^{n+1}_j&=\frac{1}{h^{n+1}_j}\Big(\int_{x_{j-\frac{1}{2}}}^{x_{j}}u(t^n,x)dx+\int_{x_{j}}^{x_{j+\frac{1}{2}}}u(t^n,x)dx\Big)\notag\\
&=\frac{1}{h^{n+1}_j}\Big(\int_{x_{j-\frac{1}{2}}}^{x_{j}}L_{j-\frac{1}{2}}(t^n,x)dx+\int_{x_{j}}^{x_{j+\frac{1}{2}}}L_{j+\frac{1}{2}}(t^n,x)dx\Big)\notag\\
&=\frac{1}{h^{n+1}_j}\Big(\frac{\Delta x}{2}(u^n_{j+\frac{1}{2}}+u^n_{j-\frac{1}{2}})+\frac{(\Delta x)^2}{8}\big(s^n_{j+\frac{1}{2}}-s^n_{j-\frac{1}{2}}\big)\Big).\label{un+1}
\end{align}
\begin{figure}[h!]
    \centering
    \begin{tikzpicture}[scale=0.8,line cap=round,line join=round,,x=1cm,y=1cm]
\draw [color=black, xstep=1cm,ystep=1cm];
\clip(-3,-5) rectangle (12,5);
\draw [line width=1.2pt,dashed] (-3,-3)-- (10,-3);
\draw [line width=1.2pt,dashed] (-1,1)-- (12,1);
\draw [line width=1.2pt,dashed] (-2,-3)-- (1,1);
\draw [line width=1.2pt,dashed] (2,-3)-- (4.45,1);
\draw [line width=1.2pt,dashed] (6,-3)-- (7.8,1);
\draw [line width=1.2pt,dashed] (10,-3)-- (11,1);
\draw [line width=2pt] (2,-1)-- (6,0);
\draw [line width=2pt] (-2,0)-- (2,-1.56);
\draw [line width=2pt] (10,-1)-- (6.016290238769366,-0.7000863220638935);
\draw [line width=1.2pt,dashed] (-2,-3)-- (-2,0);
\draw [line width=1.2pt,dashed] (2,-3)-- (2,-1);
\draw [line width=1.2pt,dashed] (6,-3)-- (6,0);
\draw [line width=1.2pt,dashed] (10.030588235294118,-2.8776470588235292)-- (10,-1);
\draw [->,line width=1.2pt,blue,dashed] (0,-3) -- (3.44,1);
\draw [->,line width=1.2pt,blue,dashed] (4,-3) -- (5,1);
\draw [->,line width=1.2pt,blue,dashed] (8,-3) -- (10,1);
\draw [line width=1.2pt] (3.32,1)-- (3.28,2.96);
\draw [line width=1.2pt] (5,1)-- (5,3);
\draw [line width=1.2pt] (10,1)-- (10,4);
\draw [line width=1.2pt] (3.28,2.96)-- (5,3);
\draw [line width=1.2pt] (5,3)-- (5,4);
\draw [line width=1.2pt] (5,4)-- (10,4);
\draw [line width=1.2pt] (3.292647793505412,2.340258118234804)-- (-1.44,2.28);
\begin{scriptsize}
\draw [fill=black] (-2,0) circle (2pt);
\draw [fill=black] (2,-1.56) circle (2pt);
\draw [fill=black] (2,-1) circle (2pt);
\draw [fill=black] (6,0) circle (2pt);
\draw [fill=black] (6,-0.7) circle (2pt);
\draw [fill=black] (10,-1) circle (2pt);
\draw (-2,-3.5) node[] {$ x^{n}_{j-\frac{3}{2}} $};
\draw (0,-3.5) node[] {$ x^{n}_{j-1} $};
\draw (2,-3.5) node[] {$ x^{n}_{j-\frac{1}{2}} $};
\draw (4,-3.5) node[] {$ x^{n}_{j} $};
\draw (6,-3.5) node[] {$ x^{n}_{j+\frac{1}{2}} $};
\draw (8,-3.5) node[] {$ x^{n}_{j+1} $};
\draw (10,-3.5) node[] {$ x^{n}_{j+\frac{3}{2}} $};
\draw (2.8,1.3) node[] {$ h^{n}_{j-1} $};
\draw (4.6,1.3) node[] {$ h^{n}_{j} $};
\draw (9.6,1.3) node[] {$ h^{n}_{j+1} $};
\draw [decorate,decoration={brace,amplitude=8pt,mirror,raise=4pt},yshift=0pt]
(-2,-2.5) -- (2,-2.5) node [black,midway,xshift=0.6cm] {\footnotesize $L_{j-1}(t^n,x)$};
\draw [decorate,decoration={brace,amplitude=8pt,mirror,raise=4pt},yshift=0pt]
(2,-2.5) -- (6,-2.5) node [black,midway,xshift=0.6cm] {\footnotesize $L_{j}(t^n,x)$};
\draw [decorate,decoration={brace,amplitude=7pt,mirror,raise=4pt},yshift=0pt]
(6,-2.5) -- (10,-2.5) node [black,midway,xshift=0.6cm] {\footnotesize $L_{j+1}(t^n,x)$};
\draw [decorate,decoration={brace,amplitude=8pt,mirror,raise=4pt},yshift=0pt]
(3.3,3.49) -- (5,3.49) node [black,midway,xshift=0.3cm] {\footnotesize $u^{n+1}_{j+\frac{1}{2}}$};
\end{scriptsize}
\end{tikzpicture}
    \caption{Geometric description of calculus of approximate $u^{n+1}_{j+\frac{1}{2}}$.}
    \label{fig:03}
\end{figure}
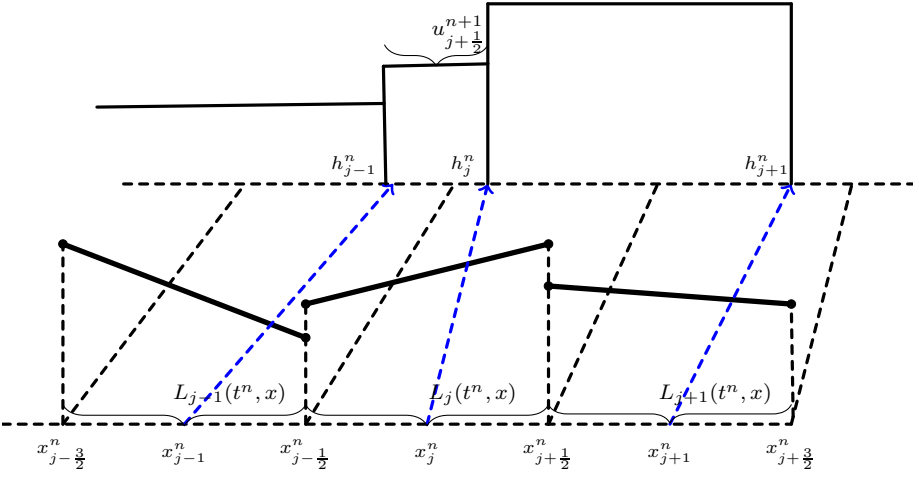
\textbf{Step 2.}
The value defined in \eqref{un+1} allows us to define the approximate value of \( u \) at time \( t^{n+1} \) as follows:
\begin{align*}
u^{n+1}_{j+\frac{1}{2}}&=\frac{1}{\Delta x}\int_{x_{j}}^{x_{j+1}}u(x,t^{n+1})dx=\frac{1}{\Delta x}\Big(\int_{x_{j}}^{\overline{x}^{n+1}_{j+\frac{1}{2}}}u^{n+1}_{j}dx+\int_{\overline{x}^{n+1}_{j+\frac{1}{2}}}^{x_{j+1}}u^{n+1}_{j+1}dx\Big)\notag\\
&=\frac{1}{\Delta x}\Big[\Big(\frac{\Delta x}{2}+f^n_{j+\frac{1}{2}}\Big)u^{n+1}_{j}+\Big(\frac{\Delta x}{2}-f^n_{j+\frac{1}{2}}\Delta t\Big)u^{n+1}_{j+1}\Big].
\end{align*} 
See Figure \ref{fig:03} for a geometric description of this step. By using the following algebraic identities:
\begin{align*}
&\frac{1}{2}\Big(h^{n+1}_{j}+\Delta t\big(f^n_{j+\frac{1}{2}}+f^n_{j-\frac{1}{2}}\big)\Big)=\frac{\Delta x}{2}+f^n_{j+\frac{1}{2}},\\
&\frac{1}{2}\Big(h^{n+1}_{j+1}-\Delta t\big(f^n_{j+\frac{3}{2}}+f^n_{j+\frac{1}{2}})\Big)=\frac{\Delta x}{2}-f^n_{j+\frac{1}{2}},
\end{align*}
we obtain
\begin{align}
u_{j+\frac{1}{2}}^{n+1}=&
\frac{1}{\Delta x}\Big[\frac{1}{2}\Big(h^{n+1}_{j}+\Delta t\big(f^n_{j+\frac{1}{2}}+f^n_{j-\frac{1}{2}}\big)\Big)u^{n+1}_{j}+\notag\\
&+\frac{1}{2}\Big(h^{n+1}_{j}-\Delta t\big(f^n_{j+\frac{3}{2}}+f^n_{j+\frac{1}{2}}\big)\Big)u^{n+1}_{j+1}\Big].\label{sche-1}
\end{align} 
Therefore, by defining the numerical flux
\begin{align}
&NF(u^n_{j-\frac{1}{2}},u^n_{j+\frac{1}{2}})=\notag\\
=&\frac{1}{4}\Big[\frac{\Delta x}{\Delta t}(u^n_{j-\frac{1}{2}}-u^n_{j+\frac{1}{2}})+\Delta x\frac{f^n_{j-\frac{1}{2}}+f^n_{j+\frac{1}{2}}}{h^{n+1}_j}(u^n_{j-\frac{1}{2}}+u^n_{j+\frac{1}{2}})\notag\\
&+\frac{(\Delta x)^2}{4}\frac{f^n_{j-\frac{1}{2}}+f^n_{j+\frac{1}{2}}}{h^{n+1}_j}\Big(s^n_{j-\frac{1}{2}}-s^n_{j+\frac{1}{2}}\Big)+\frac{(\Delta x)^2}{4\Delta t}\Big(s^n_{j-\frac{1}{2}}+s^n_{j+\frac{1}{2}}\Big)\Big],\label{Flux-n1}
\end{align}
and using this function in \eqref{sche-1}, 
we obtain the following scheme in conservative form
\begin{equation}\label{numerical}
u^{n+1}_{j+\frac{1}{2}}=u^n_{j+\frac{1}{2}}-\frac{\Delta t}{\Delta x}\Big(NF(u^n_{j+\frac{1}{2}},u^n_{j+\frac{3}{2}})-NF(u^n_{j-\frac{1}{2}},u^n_{j+\frac{1}{2}})\Big). 
\end{equation}
\begin{rmk}
By using the condition provided in \eqref{CFL-0}, the values \( \frac{\Delta x}{\Delta t} \) can be removed from \eqref{Flux-n1} by introducing the parameters:
\begin{equation}\label{bj}
b^n_{j}=b_{j}\Big(f^n_{j-\frac{1}{2}},f^n_{j+\frac{1}{2}}\Big). 
\end{equation}
In fact, notice that the condition \eqref{CFL-0}, in terms of \eqref{f_j}, can be written as
\begin{equation}\label{cfl-0}
-\frac{\Delta x}{\Delta t} \leq f^n_{j-\frac{1}{2}} + f^n_{j+\frac{1}{2}} \leq \frac{\Delta x}{\Delta t},
\end{equation}
therefore, there exists a value \( b = b\big(f^n_{j-\frac{1}{2}}, f^n_{j+\frac{1}{2}}\big) \in (0,1) \), for any choice of \( \Delta t \), such that
\begin{equation}\label{sum-cf}
\frac{\Delta x}{\Delta t} = \frac{f^n_{j-\frac{1}{2}} + f^n_{j+\frac{1}{2}}}{2b-1} = b^n_{j}.
\end{equation}
Finally, by construction, we obtain the inequality 
\begin{equation}\label{b-ineq}
|f^n_{j-\frac{1}{2}} + f^n_{j+\frac{1}{2}}| \leq b^n_{j}.
\end{equation}
\end{rmk}

\textbf{Step 4.} Replacing \eqref{bj} into \eqref{Flux-n1}, we obtain
\begin{align}
&NF(u^n_{j-\frac{1}{2}}, u^n_{j+\frac{1}{2}}) = \notag \\
&\frac{1}{4} \Big[b^n_{j}(u^n_{j-\frac{1}{2}} - u^n_{j+\frac{1}{2}}) 
+ \Delta x \frac{f^n_{j-\frac{1}{2}} + f^n_{j+\frac{1}{2}}}{h^{n+1}_j}(u^n_{j-\frac{1}{2}} + u^n_{j+\frac{1}{2}}) \notag\\
&\quad + \frac{(\Delta x)^2}{4} \frac{f^n_{j-\frac{1}{2}} + f^n_{j+\frac{1}{2}}}{h^{n+1}_j} \Big(s^n_{j-\frac{1}{2}} - s^n_{j+\frac{1}{2}}\Big) 
+ \frac{\Delta x}{4} b^n_{j} \Big(s^n_{j-\frac{1}{2}} + s^n_{j+\frac{1}{2}}\Big)\Big] \notag\\
&= \frac{1}{4} \Big\{b^n_{j} \Big(\big(u^n_{j-\frac{1}{2}} + \frac{\Delta x}{4}s^n_{j-\frac{1}{2}}\big) 
- \big(u^n_{j+\frac{1}{2}} - \frac{\Delta x}{4}s^n_{j+\frac{1}{2}}\big)\Big) \notag\\
&\quad + \Delta x \frac{f^n_{j-\frac{1}{2}} + f^n_{j+\frac{1}{2}}}{h^{n+1}_j} 
\Big(\big(u^n_{j-\frac{1}{2}} + \frac{\Delta x}{4}s^n_{j-\frac{1}{2}}\big) 
+ \big(u^n_{j+\frac{1}{2}} - \frac{\Delta x}{4}s^n_{j+\frac{1}{2}}\big)\Big)\Big\}.
\label{Flux-1-n2-0}
\end{align}
Let 
\begin{equation}\label{pm-eq1}
\begin{aligned}
& u^{-,n}_{j}=u^n_{j-\frac{1}{2}}+\frac{\Delta x}{4}s^n_{j-\frac{1}{2}},\\
& u^{+,n}_{j}=u^n_{j+\frac{1}{2}}-\frac{\Delta x}{4}s^n_{j+\frac{1}{2}},
\end{aligned}
\end{equation}
we can observe that these values are related to the jump of the function \( L(x) \) 
defined in \eqref{L-1} (see Figure \ref{fig:4} for a geometric description). 
Replacing \eqref{pm-eq1} into \eqref{Flux-1-n2-0}, 
we obtain the following representation of the flux \( NF \):
\begin{align*}
&NF(u^n_{j-\frac{1}{2}},u^n_{j+\frac{1}{2}})=NF(u^{-,n}_j,u^{+,n}_j)\\
&=\frac{1}{4}\Big\{b^n_{j}\Big(u^{-,n}_{j}-u^{+,n}_{j+1}\Big)+\Delta x\frac{f^n_{j-\frac{1}{2}}+f^n_{j+\frac{1}{2}}}{h^{n+1}_{j}}\Big(u^{-,n}_{j}+u^{+,n}_{j+1}\Big)\Big\}.
\notag
\end{align*}

Since \( \lim_{\Delta t \to 0} h^{n+1}_{j} = \Delta x \), then if we take the limit \( \Delta t \to 0 \) in \eqref{numerical}, we obtain the following semidiscrete scheme:
\begin{equation}\label{Semi}
\frac{d}{dt}u_{j+\frac{1}{2}}(t)=-\frac{1}{\Delta x}\Big(\mathcal{F}(u^{-}_{j+1}(t),u^{+}_{j+1}(t))-\mathcal{F}(u^-_{j}(t),u^{+}_{j}(t))\Big), 
\end{equation}
where
\begin{align}
\mathcal{F}(u^-_{j}(t), u^+_{j}(t)) = 
&\frac{1}{4} \Big[ b_{j}(t) \Big(u^-_{j}(t) - u^+_{j}(t)\Big) \notag\\
&\quad + \big(f_{j-\frac{1}{2}}(t) + f_{j+\frac{1}{2}}(t)\big) \Big(u^-_{j}(t) + u^+_{j}(t)\Big) \Big],\label{num_flux}
\end{align}

and
\begin{equation}
u^{-}_{j}(t)=u_{j-\frac{1}{2}}(t)+\frac{\Delta x}{4}s_{j-\frac{1}{2}}(t), \ u^{+}_{j}(t)=u_{j+\frac{1}{2}}(t)-\frac{\Delta x}{4}s_{j+\frac{1}{2}}(t).\label{upm} 
\end{equation}
\begin{figure}
    \centering
    \begin{tikzpicture}[line cap=round,line join=round]
\draw [color=black,dash pattern=on 1pt off 1pt, xstep=1.0cm,ystep=1.0cm];
\clip(-3,-4.02) rectangle (6.98,2.94);
\draw (-3,-3)-- (7,-3);
\draw (-2,-1.66)-- (2,-0.46);
\draw (2,0.2)-- (6,1.5);
\draw (-2.24,-3) node[anchor=north west] {$ x_{j-1} $};
\draw (-0.28,-3) node[anchor=north west] {$ x_{j-\frac{1}{2}} $};
\draw (1.86,-3) node[anchor=north west] {$ x_j $};
\draw (3.78,-3) node[anchor=north west] {$ x_{j+\frac{1}{2}} $};
\draw (5.76,-3) node[anchor=north west] {$ x_{j+1} $};
\draw (-0.7,-0.22) node[anchor=north west] {$ u^n_{j-\frac{1}{2}} $};
\draw (3.22,1.78) node[anchor=north west] {$u^n_{j+\frac{1}{2}} $};
\draw (2.16,-0.2) node[anchor=north west] {$ L^+_{j-\frac{1}{2}}= u_{j-\frac{1}{2}}(t)+\frac{\Delta x}{2}d_{j-\frac{1}{2}} $};
\draw (2.16,0.5) node[anchor=north west] {$ L^-_{j+\frac{1}{2}}= u_{j+\frac{1}{2}}(t)-\frac{\Delta x}{2}d_{j+\frac{1}{2}} $};
\begin{scriptsize}
\draw [dashed ](2.0,-3)--(2.0,3);
\fill [color=black] (-2,-3) circle (1.5pt);
\fill [color=black] (2,-3) circle (1.5pt);
\fill [color=black] (6,-3) circle (1.5pt);
\fill [color=blue] (0,-1) circle (1.7pt);
\fill [color=red] (2,-0.46) circle (1.7pt);
\fill [color=red] (2,0.2) circle (1.7pt);
\fill [color=blue] (4.02,0.86) circle (1.7pt);
\end{scriptsize}
\end{tikzpicture}
    \caption{Graphical description of the 
jumps of the function $L(x)$ at $x=x_j$.}
    \label{fig:4}
\end{figure}
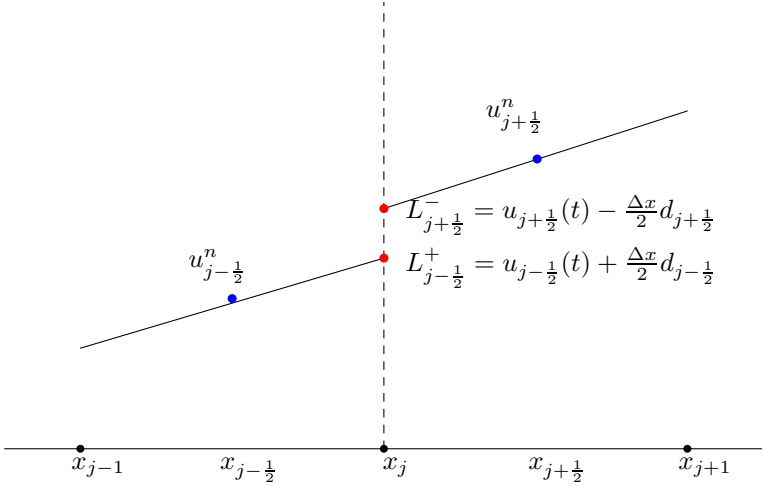

\begin{rmk}\label{Rem2.38}
\begin{enumerate}
\item[1)] By using \eqref{upm}, we obtain the 
following representation of \( u_{j+\frac{1}{2}}(t) \):
\begin{equation}\label{Iden-u}
u_{j+\frac{1}{2}}(t) = \frac{1}{2}\Big(u^-_{j}(t) + u^{+}_{j+1}(t)\Big).
\end{equation}

\item[2)] If we define the functions
\begin{align}
P_{j}(t) &= \frac{1}{4}\Big(f_{j-\frac{1}{2}}(t) + f_{j+\frac{1}{2}}(t) + b_{j}(t)\Big), \label{Pj}\\
N_{j}(t) &= \frac{1}{4}\Big(f_{j-\frac{1}{2}}(t) + f_{j+\frac{1}{2}}(t) - b_{j}(t)\Big), \label{Nj}
\end{align}
then, from conditions \eqref{cfl-0} and \eqref{sum-cf}, 
it follows that \( P_{j+\frac{1}{2}}(t) \geq 0 \) and 
\( N_{j+\frac{1}{2}}(t) \leq 0 \), and we have the 
following equivalent formulation for the ODE \eqref{Semi}:
\begin{align}\label{Semi-eq1}
\frac{d}{dt}u_{j+\frac{1}{2}}(t) &= \frac{1}{\Delta x} 
\Big(u^{-}_{j}(t)P_{j}(t) 
- \big(u^{-}_{j+1}(t)P_{j+1}(t) - u^{+}_{j}(t)N_{j}(t)\big) 
- u^+_{j+1}(t)N_{j+1}(t)\Big). 
\end{align}
\end{enumerate}
\end{rmk}

\begin{rmk}\label{Rem2.39}
By introducing
\begin{equation}\label{G-func}
\mathcal{G}_j(u,v) = \frac{1}{4}\Big(b_{j}(u-v) + \big(f_{j-\frac{1}{2}}(t) + f_{j+\frac{1}{2}}(t)\big)(u+v)\Big),
\end{equation}
\eqref{Semi-eq1} can be written compactly as
\begin{equation}\label{Semi-pm}
\frac{d}{dt}u_{j+\frac{1}{2}}(t) = -\frac{1}{\Delta x}\Big(\mathcal{G}_{j+1}(u^{-}_{j+1}(t), u^{+}_{j+1}(t)) - \mathcal{G}_j(u^{-}_{j}(t), u^{+}_{j}(t))\Big). 
\end{equation}
Moreover, the function \( \mathcal{G}_j(u,v) \) defined in \eqref{G-func} satisfies the following properties:
\begin{itemize}
\item[\textbf{G$_1$.}] \textit{Consistency}: In line with classical approximations for scalar conservation laws (see, e.g., \cite{MCAM80}), 
for \( \mathcal{G} \) we have
\begin{equation*}
\mathcal{G}_{j}(u,u) =
\begin{cases}
k_Lf(u), & \text{if}\ j<0,\\
\frac{1}{2}(k_L + k_R)f(u), & \text{if}\ j=0,\\
k_Rf(u), & \text{if}\ j>0.\\
\end{cases}
\end{equation*}
\item[\textbf{G$_2$.}] \textit{Non-decreasing in the first argument}: Indeed, by the definition of \( b_{j}(t) \), we have
\[
\partial_u \mathcal{G}_j(u,v) = \frac{1}{4}\Big(b_{j+\frac{1}{2}} + \big(f_{j-\frac{1}{2}}(t) + f_{j+\frac{1}{2}}(t)\big)\Big) \geq 0.
\]
\item[\textbf{G$_3$.}] \textit{Non-increasing in the second argument}: By \eqref{b-ineq}, it follows that
\[
\partial_v \mathcal{G}_j(u,v) = \frac{1}{4}\Big(-b_{j}(t) + \big(f_{j-\frac{1}{2}}(t) + f_{j+\frac{1}{2}}(t)\big)\Big) \leq 0.
\]
\end{itemize}
\end{rmk}
\section{Convergence analysis of the semidiscrete Lagrangian-Eulerian}\label{aansd}
In this section, we prove the convergence of the sequence of approximate solutions obtained by \eqref{Semi} and show that the limit function satisfies \eqref{3-loc-weak}. 
The strategy involves establishing uniform bounds for the \( L^1 \) and 
\( L^{\infty} \) norms, providing a suitable modulus of continuity, 
and ensuring \( L^1 \)-continuity in time. These properties enable us 
to apply \( L^1 \)-compactness tools and demonstrate the existence of a 
convergent subsequence. Considering the spatial discretization 
introduced at the beginning of Subsection \ref{Subsec-2.1}, 
we define the following step function
\begin{equation}\label{Fin-app}
u^{\Delta x}(t,x)=\sum_{j}u_{j+\frac{1}{2}}(t)\mathbbm{1}_{I_{j}}(x), \quad I_{j}=(x_{j},x_{j+1}).
\end{equation}
We also consider the piecewise constant approximation of the initial data given by
\begin{equation*}
u^{\Delta x}_0(x)=\sum_{j}u_{j+\frac{1}{2}}(0)\mathbbm{1}_{I_{j}}(x), \quad I_{j}=(x_{j},x_{j+1}),
\end{equation*}
where $u_{j+\frac{1}{2}}(0)=\frac{1}{\Delta x}\int_{I_j}u_0(x)dx$. 
The approximate solution for \eqref{1-1}-\eqref{1-1-D}
can be constructed by solving the infinite system of ordinary differential 
equations \eqref{Semi} written in compact  form as
\begin{align}
& \frac{d}{dt}u^{\Delta }(t)=G(u^{\Delta}(t)),\label{Cauchy-1}\\
&u^{\Delta x}(0)=u^{\Delta}_{0}.\label{Cauchy-Data}
\end{align} 
where
\begin{equation}\label{Flux-step}
G(u^{\Delta}(t,x))_j=-\frac{1}{\Delta x}\Big(\mathcal{G}_{j+1}(u^{-}_{j+1}(t),u^{+}_{j+1}(t))-\mathcal{G}_{j}(u^-_{j}(t),u^{+}_{j}(t))\Big).
\end{equation}
In the next results, we will use the following discretized versions of the \( L^1 \) norm and Total Variation:
\begin{align}
& \|u^{\Delta x}(t,\cdot)\|_{1} = \Delta x \sum_{j} |u_{j+\frac{1}{2}}(t)|, \label{L1-norm}\\
& TV(|u^{\Delta x}(t,\cdot)|) = \sum_{j} |u_{j+\frac{1}{2}}(t) - u_{j-\frac{1}{2}}(t)|. \label{TV-norm}
\end{align}
\begin{thm}\label{Thm-exis}
Let \( \Delta x > 0 \) be a spatial discretization parameter. 
Suppose that \( u^{\Delta x}_0 \in L^1(\mathbb{R}) \cap BV(\mathbb{R}) \). 
Then, there exists a time horizon \( T > 0 \) such that the initial value 
problem \eqref{Cauchy-1}-\eqref{Cauchy-Data} has a unique solution in \( C^{1}([0,T]; L^{1}(\mathbb{R})) \).
\end{thm}
\begin{proof}
Notice that in the definition of the function \( G \), given in \eqref{Flux-step}, 
the operators \( u^{\pm}_{j+1}(t) \) depend on the values 
\( u_{j+\frac{1}{2}}(t) \). Therefore, by using \eqref{upm}, we have that
\begin{align}
&|u^{\pm}_{j+1}(t)-v^{\pm}_{j+1}(t)|\notag\\
&\leq |u_{j+\frac{1}{2}}(t)-v_{j+\frac{1}{2}}(t)|+\frac{1}{4}\Big(|u_{j+\frac{1}{2}}(t)-u_{j-\frac{1}{2}}(t)|+|v_{j+\frac{1}{2}}(t)-v_{j-\frac{1}{2}}(t)|\Big)\label{Exis}.
\end{align}
Moreover, the functions \( \mathcal{G}_{j}(u,v) \) defined in \eqref{G-func} are Lipschitz in all variables. Hence, by rewriting
\begin{align}
&\mathcal{G}_{j+1}(u^{-}_{j+1}(t),u^{+}_{j+1}(t))-\mathcal{G}_{j+1}(v^{-}_{j+1}(t),v^{+}_{j+1}(t))\notag\\
&=\mathcal{G}_{j+1}(u^{-}_{j+1}(t),u^{+}_{j+1}(t))-\mathcal{G}_{j+1}(v^{-}_{j+1}(t),u^{+}_{j+1}(t))\notag\\
&+\mathcal{G}_{j+1}(v^{-}_{j+1}(t),u^{+}_{j+1}(t))-\mathcal{G}_{j+1}(v^{-}_{j+1}(t),v^{+}_{j+1}(t)).
\label{Exis-2}
\end{align}
and taking absolute value in \eqref{Exis-2} we obtain
\begin{align}
&|\mathcal{G}_{j+1}(u^{-}_{j+1}(t),u^{+}_{j+1}(t))-\mathcal{G}_{j+1}(v^{-}_{j+1}(t),v^{+}_{j+1}(t))|\notag\\
&\leq |\mathcal{G}_{j+1}(u^{-}_{j+1}(t),u^{+}_{j+1}(t))-\mathcal{G}_{j+1}(v^{-}_{j+1}(t),u^{+}_{j+1}(t))|\notag\\
&+|\mathcal{G}_{j+1}(v^{-}_{j+1}(t),u^{+}_{j+1}(t))-\mathcal{G}_{j+1}(v^{-}_{j+1}(t),v^{+}_{j+1}(t))|\notag\\
&\leq \text{Lip}_u(\mathcal{G})|u^{-}_{j+1}(t)-v^{-}_{j+1}(t)|+\text{Lip}_v(\mathcal{G})|u^{+}_{j+1}(t)-v^{+}_{j+1}(t)|,
\label{Exis-3}\end{align}
where \( \text{Lip}_u(\mathcal{G}) \) and \( \text{Lip}_v(\mathcal{G}) \) 
represent the Lipschitz constants on the first and the 
second variables, respectively. Now, consider the integral operator:
\begin{equation}\label{I-operator}
\Psi(u^{\Delta x})(t,x)=u^{\Delta x}_0(x)+\int_{0}^tG(u^{\Delta x}(s,x))ds.
\end{equation}
By \eqref{Exis-2} and \eqref{Exis-3}, the flux function \( G \) defined in 
\eqref{Flux-step} is a Lipschitz function. Therefore, by \eqref{Exis}, 
the integral operator \eqref{I-operator} is Lipschitz in the \( L^1 \)-norm. 
By an application of the Picard-Lindelöf Theorem \cite[Theorem 3.A]{Zeidler}, 
we obtain the desired result.
\end{proof}
\subsection{$L^{\infty}$-estimates}\label{subsec3-1}
In this section, we provide uniform bounds for the family \( \big\{u_{j+\frac{1}{2}}(t)\big\}_{j\in \mathbb{Z}} \) obtained by the semi-discrete Lagrangian-Eulerian scheme \eqref{Semi}. 
To achieve this result, we use the equivalent ODE \eqref{Semi-pm}.

\begin{prop}\label{p-maximum}
Let \( \big\{u_{j+\frac{1}{2}}(t)\big\}_{j\in \mathbb{Z}} \) be the family of  
solutions to the ODE \eqref{Semi-pm}. If the initial data \( u_0 \) 
lies within the interval \( I = [0,1] \), then 
\begin{equation}\label{Max-bounds}
u_{j+\frac{1}{2}}(t) \in [0,1], \quad \text{for all}\ j \in \mathbb{Z},\ \text{and}\ t \in [0,T].
\end{equation}
\end{prop}

\begin{proof}
For \( \Delta t > 0 \) small enough, by using \eqref{Iden-u}, we 
can write the ODE \eqref{Semi-pm} in the following form:
\begin{align}\label{Max-1}
&u_{j+\frac{1}{2}}(t+dt)=\frac{1}{2}\Big(u^{-}_{j+1}(t)+u^{+}_{j}(t)\Big)\notag\\
&-\frac{dt}{\Delta x}\Big(\mathcal{G}_{j+1}(u^{-}_{j+1}(t),u^{+}_{j+1}(t))-\mathcal{G}_j(u^{-}_{j}(t),u^{+}_{j}(t))\Big)+ O(dt)\notag\\
&=\frac{1}{2}u^{-}_{j+1}(t)-\frac{dt}{\Delta x}\mathcal{G}_{j+1}(u^{-}_{j+1}(t),u^{+}_{j+1}(t))\notag\\
&+\frac{1}{2}u^{+}_{j}(t)+\frac{dt}{\Delta x}\mathcal{G}_j(u^{-}_{j}(t),u^{+}_{j}(t))+ O(dt).
\end{align}
Let
\[
R(u,v)=\frac{1}{2}u-\frac{dt}{\Delta x}\mathcal{G}_{j+1}(u,v),
\]
and
\[
S(u,v)=\frac{1}{2}v+\frac{dt}{\Delta x}\mathcal{G}_j(u,v),
\]
by using the properties $G_2$ and $G_3$ from Remark \eqref{Rem2.39}, 
we have that
\[
\frac{\partial}{\partial u}R(u,v)=\frac{1}{2}-\frac{dt}{\Delta x}\frac{\partial}{\partial u}\mathcal{G}_{j+1}(u,v)\geq 0,
\]
and
\[
\frac{\partial}{\partial v}S(u,v)=\frac{1}{2}+\frac{dt}{\Delta x}\frac{\partial}{\partial v}\mathcal{G}_{j}(u,v)\geq 0.
\]
Finally, let
\[
H(u_{j-\frac{1}{2}},u_{j+\frac{1}{2}},u_{j+\frac{3}{2}})=R(u_{j+\frac{3}{2}}(t),u_{j+\frac{1}{2}}(t))+S(u_{j+\frac{1}{2}}(t),u_{j-\frac{1}{2}}(t)),\]
we can rewrite \eqref{Max-1} as
\[
u_{j+\frac{1}{2}}(t+dt)=H(u_{j-\frac{1}{2}},u_{j+\frac{1}{2}},u_{j+\frac{3}{2}})+ O(dt),
\]
where the flux satisfies
\[
\frac{\partial H(u_{j-\frac{1}{2}},u_{j+\frac{1}{2}},u_{j+\frac{3}{2}})}{\partial u_{j+\frac{1}{2}}} \geq 0
\]
and conclude that $u_{j+\frac{1}{2}}(t+dt)\in [0,1]$. 
Finally, considering the partition $t_{k}=k dt$, an induction  
argument applied to the initial value problem \eqref{Semi-pm} with initial datum
\[
u_{j+\frac{1}{2}}(t_{k})=u_{j+\frac{1}{2}}(t_{k-1}),
\]
produces the desired  result.
\end{proof}
\subsection{$L^1$-estimates}\label{subsec-3.2}
To provide \( L^1 \)-estimates for the family \( \{u_{j+\frac{1}{2}}(t)\}_{j\in \mathbb{Z}} \), 
it is computationally more convenient to work with the ODE \eqref{Semi-eq1}.
\begin{prop}\label{L1-est}
Let \( u_0 \in L^{1}(\mathbb{R}) \) be the initial data, and 
\( \{u_{j+\frac{1}{2}}(t)\}_{j} \) the corresponding solutions  
computed by the semidiscrete Lagrangian-Eulerian scheme \eqref{Semi}. Then we have
\begin{equation}\label{L-1bound}
\sum_{j}|u_{j+\frac{1}{2}}(t)| \leq \sum_{j}|u_{j+\frac{1}{2}}(0)|, 
\end{equation}
for all \( t \in (0,T] \).
\end{prop}
\begin{proof}
Recalling Remark \eqref{Rem2.38}, we consider the equivalent ODE 
proposed in \eqref{Semi-eq1}. By Theorem \ref{Thm-exis}, the function 
\( u^{\Delta x}(t,x) \in C^{1}([0,T]; L^{1}(\mathbb{R})) \). Then, for 
any sufficiently small \( \Delta t > 0 \), we can apply a first-order 
Taylor approximation and obtain
\begin{align}
& u_{j+\frac{1}{2}}(t+dt)=u_{j+\frac{1}{2}}(t)+\frac{dt}{\Delta x}\Big(u^{-}_{j}(t)P_{j-\frac{1}{2}}(t)
-\big(u^{-}_{j+1}(t)P_{j+\frac{1}{2}}(t)-\notag\\&-u^{+}_{j}(t)N_{j-\frac{1}{2}}(t)\big)-u^+_{j+1}(t)N_{j+\frac{1}{2}}(t)\Big)+ dt\cdot r_{j+\frac{1}{2}}(t,dt) \label{L1-1}
\end{align}
By using the identity \eqref{upm}, we write the ODE \eqref{L1-1} 
in the following  convenient  form
\begin{align}
& u_{j+\frac{1}{2}}(t+dt)= u^{-}_{j+1}(t)\Big(\frac{1}{2}-\frac{dt}{\Delta x}P_{j+1}(t)\Big)+u^{+}_{j}(t)\Big(\frac{1}{2}+\frac{dt}{\Delta x}N_{j}(t)\Big)\notag\\
&+\frac{dt}{\Delta x}u^-_{j}(t)P_j(t)-\frac{dt}{\Delta x}u^+_{j+1}(t)N_{j+1}(t)+ dt\cdot r_{j+\frac{1}{2}}(t,dt). \label{L1-2}
\end{align}
From Remark \ref{Rem2.38}, the functions \( P_{j}(t) \) and \( N_{j}(t) \), 
defined by \eqref{Pj} and \eqref{Nj}, respectively, satisfy
\[
P_{j}(t) \geq 0 \quad \text{and} \quad N_{j}(t) \leq 0.
\]
Moreover, these depend on 
the values \( u_{j+\frac{1}{2}}(t) \). Therefore, by \eqref{Max-bounds} 
from Proposition \ref{p-maximum}, they are uniformly 
bounded. Then, we can choose a sufficiently small \( \Delta t > 0 \) such that the following 
conditions hold:
\begin{equation}\label{L1-cond}
\frac{1}{2} - \frac{\Delta t}{\Delta x} P_{j+1}(t) < 0, \quad \text{and} \quad \frac{1}{2} + \frac{\Delta t}{\Delta x} N_{j+1}(t) > 0.
\end{equation}
Taking the absolute value in \eqref{L1-2}, and using the condition \eqref{L1-cond},  we obtain
\begin{align}
& |u_{j+\frac{1}{2}}(t+dt)|\leq |u^{-}_{j+1}(t)|\Big(\frac{1}{2}-\frac{dt}{\Delta x}P_{j+1}(t)\Big)+|u^{+}_{j}(t)|\Big(\frac{1}{2}+\frac{dt}{\Delta x}N_{j}(t)\Big)\notag \\
&+\frac{dt}{\Delta x}|u^-_{j}(t)|P_j(t)-\frac{dt}{\Delta x}|u^+_{j+1}(t)|N_{j+1}(t)+dt\cdot r_{j+\frac{1}{2}}(t,dt).\label{L1-3}
\end{align}
Summing over $j$ in equation\eqref{L1-3}, we get
\begin{align}
& \sum_{j}| u_{j+\frac{1}{2}}(t+dt)|\leq \frac{1}{2}\sum_{j}\big(|u^{-}_{j+1}(t)|+|u^{+}_{j}(t)|\big)\label{L1-4}\\
&+ \frac{dt}{\Delta x}\Big(-\sum_{j}|u^{-}_{j+1}(t)|P_{j+1}(t)+\sum_{j}|u^-_{j}(t)|P_j(t)\Big)\label{L1-4-1}\\
&\frac{dt}{\Delta x}\Big(\sum_{j}|u^{+}_{j}(t)|N_{j}(t)-\sum_{j}|u^+_{j+1}(t)|N_{j+1}(t)\Big)\label{L1-4-2}\\
&+dt\sum_{j} r_{j+\frac{1}{2}}(t,dt).\notag
\end{align}
Notice that the terms \eqref{L1-4-1} and \eqref{L1-4-2} vanish, 
then \eqref{L1-4} reduces to  
\begin{equation}\label{L1-5}
\sum_{j}| \frac{d}{dt}{u}_{j+\frac{1}{2}}(t)|\leq \sum_{j}\big(|u^{-}_{j+1}(t)|+|u^{+}_{j}(t)|\big)+dt\sum_{j} r_{j+\frac{1}{2}}(t,dt)
\end{equation}
Now, we decompose the sum on the right hand side of \eqref{L1-5} in the following sets
\begin{enumerate}[1).]
\item The set $A_0$ where $u_{j+\frac{1}{2}}(t)=0$, in this  case we have
\[
|u^{-}_{j+1}(t)|=|u^{+}_{j}(t)|.
\]
\item The set $A^{+}$ where $u_{j+\frac{1}{2}}(t)>0$,
\item The set $A^{-}$ where $u_{j+\frac{1}{2}}(t)>0$.
\end{enumerate}
Then we can choose $\Delta x$ small enough so that
(by \eqref{upm}), we have
\begin{align}
|u^{+}_{j}(t)|&=|u_{j+\frac{1}{2}}(t)|-\frac{\Delta x}{4}\text{sgn}(u_{j+\frac{1}{2}}(t))s_{j+\frac{1}{2}}(t),\label{L1-6-1}\\
|u^{-}_{j+1}(t)|&=|u_{j+\frac{1}{2}}(t)|+\frac{\Delta x}{4}\text{sgn}(u_{j+\frac{1}{2}}(t))s_{j+\frac{1}{2}}(t),\label{L1-6-2} 
\end{align}
for values in the sets \( A^{+} \) or \( A^{-} \).
Replacing \eqref{L1-6-1} and \eqref{L1-6-2} into \eqref{L1-5}, we obtain
\begin{align}
&\sum_{j}|u_{j+\frac{1}{2}}(t+dt)|\leq \sum_{j}|u_{j+\frac{1}{2}}(t)|+dt\|r(t,dt)\|_{\ell^1}.\label{L1-7}
\end{align}
Now we consider the partition \( t_k = k \cdot \Delta t \), \( t_N = N \cdot \Delta t = t \). 
Applying \eqref{L1-7} recursively, we obtain
\begin{align}
&\sum_{j}|u_{j+\frac{1}{2}}(t)|\leq \sum_{j}|u_{j+\frac{1}{2}}(0)|+dt\sum_{k}\|r(t_k,dt)\|_{\ell^1}\notag\\
&\leq \sum_{j}|u_{j+\frac{1}{2}}(0)|+ N\cdot dt \max_{k}\big(\|r(t_k,dt)\|_{\ell^1}\big)\notag\\
&\leq \sum_{j}|u_{j+\frac{1}{2}}(0)|+ T \max_{k}\big(\|r(t_k,dt)\|_{\ell^1}\big)\label{L1-8}
\end{align}
Taking limit  $dt\to 0$ in \eqref{L1-8}, finally we get
\begin{align*}
&\sum_{j}|u_{j+\frac{1}{2}}(t)|\leq \sum_{j}|u_{j+\frac{1}{2}}(0)|
\end{align*}
\end{proof}

\subsection{Non-uniform TV-estimates} 

In the Introduction \ref{Intro}, 
we provided an example in which solutions of scalar conservation 
laws with discontinuous flux can increase the total variation. 
In this subsection, we provide an estimate of the increase in total variation for an approximate solution computed through \eqref{Semi-pm}.

\medskip
\begin{prop}\label{TV-nounif} Let $u_{0}\in BV(\mathbb{R})$ and
 $\{u_{j+\frac{1}{2}}(t)\}_{j\in \mathbb{Z}}$ be the family 
of solutions to the ODE \eqref{Semi-pm}.  
Then there exists a continuous function  $B$ such that
\begin{equation}\label{TV-b}
TV(u^{\Delta x}(t))\leq TV(u_0)+\frac{1}{\Delta x}\int_{0}^{t}B(s)ds .
\end{equation}
\end{prop}
\begin{proof}
Notice that, due to the definition of the values \( f_{j+\frac{1}{2}}(t) \) 
given in \eqref{f_j}, the flux function \( \mathcal{G}_{j} \) defined in 
\eqref{G-func} exhibits different behavior depending on the values of 
the index \( j \in \mathbb{Z} \). To account for this
dependence, we consider the following steps.\\[0.2in]

\textbf{Step 1.} When \( j > 1 \) or \( j < -1 \), the functions \( f_{j+\frac{1}{2}}(t) \) are given by:
\[
f_{j+\frac{1}{2}}(t) =
\begin{cases}
    k_R(1 - u_{j+\frac{1}{2}}(t)), & \text{if } j > 1,\\
    k_L(1 - u_{j+\frac{1}{2}}(t)), & \text{if } j < -1.
\end{cases}
\]
Recalling \eqref{Semi-pm}, we have
\begin{equation}\label{TV-1-0}
\frac{d}{dt}u_{j+\frac{1}{2}}(t)=-\frac{1}{\Delta x}\Big(\mathcal{G}_{j+1}(u^-_{j+1}(t),u^+_{j+1}(t))-\mathcal{G}_j(u^-_{j}(t),u^+_{j}(t))\Big).
\end{equation}
We rewrite the difference of fluxes in \eqref{TV-1-0} as follows:
\begin{align}
&\mathcal{G}_{j+1}(u^-_{j+1}(t),u^+_{j+1}(t))-\mathcal{G}_{j+1}(u^-_{j}(t),u^+_{j}(t))\label{Tv2-1}\\
&=\Big(\mathcal{G}_{j+1}(u^-_{j+1}(t),u^+_{j+1}(t))\notag\\
&-\mathcal{G}_{j+1}(u^-_{j}(t),u^+_{j+1}(t))\Big)-\Big(\mathcal{G}_{j+1}(u^-_{j}(t),u^+_{j+1}(t))-\mathcal{G}_{j+1}(u^-_{j}(t),u^-_{j}(t))\Big).\notag
\end{align}
Then we set
\begin{equation*}
\Delta u_{j+1}(t)=u_{j+\frac{3}{2}}(t)-u_{j+\frac{1}{2}}(t),
\end{equation*}
and
\begin{align}
& A_{j+1}(t)=\frac{1}{\Delta u_{j+1}(t)}\Big(\mathcal{G}_{j+1}(u^-_{j+1}(t),u^+_{j+1}(t))-\mathcal{G}_{j+1}(u^-_{j+1}(t),u^+_{j}(t))\Big)\label{A},\\
&B_{j}(t)=\frac{1}{\Delta u_{j}(t)}\Big(\mathcal{G}_{j}(u^-_{j+1}(t),u^+_{j}(t))-\mathcal{G}_{j}(u^-_{j}(t),u^+_{j}(t))\Big).\label{B}
\end{align}
Replacing these values into \eqref{Tv2-1}, we can write \eqref{TV-1-0} in the  following form
\begin{equation}\label{Semi-TV}
\frac{d}{dt}u_{j+\frac{1}{2}}(t)=-\frac{1}{\Delta x}A_{j+1}(t)\Delta u_{j+1}(t)+\frac{1}{\Delta x}B_{j}(t)\Delta u_{j}(t),
\end{equation}
We use \eqref{Semi-TV} to compute $\frac{d}{dt}\Delta u_{j+1}(t)$  and obtain
\begin{align}
&\frac{d}{dt}(\Delta u_{j+1}(t))=\label{diff}\\
&-\frac{1}{\Delta x}\Big(A_{j+2}(t)\Delta u_{j+2}(t)-B_{j+1}(t)\Delta u_{j+1}(t)\Big)\notag\\
&+\frac{1}{\Delta x}\Big(A_{j+1}(t)\Delta u_{j+1}(t)-B_{j}(t)\Delta u_{j}(t)\Big).\notag
\end{align}
We introduce
\[
\sigma_{j+1}(t)=\text{sgn}\big(\Delta u_{j+1}(t)\big),
\]
 then multiply by $\sigma_{j+1}(t)$ in \eqref{diff} and sum by parts to obtain
\begin{align}
&\sum_{j}\frac{d}{dt}|\Delta u_{j+1}(t)|\label{diff-2}\\
&=-\frac{1}{\Delta x}\sum_{j}\Big\{\Big(\sigma_{j}(t)-\sigma_{j+1}(t)\Big)A_{j+1}(t)\Delta u_{j+1}(t)\notag\\
&+\Big(\sigma_{j+1}(t)-\sigma_{j}(t)\Big)B_{j}(t)\Delta u_{j}(t)\Big\}\notag\\
&=\frac{1}{\Delta x}\sum_{j}\Big\{\Big(1-\sigma_{j}(t)\cdot\sigma_{j+1}(t)\Big)A_{j+1}(t)|\Delta u_{j+1}(t)|\notag\\
&+\Big(1-\sigma_{j+1}(t)\cdot\sigma_{j}(t)\Big)B_{j}(t)|\Delta u_{j}(t)|\Big\}.\notag
\end{align}
Now notice that \( \sigma_{j+1}(t) \cdot \sigma_{j}(t) = \pm 1 \). 
Then we can divide the sum \eqref{diff-2} into two sets, namely:
\begin{enumerate}[1).]
\item The set of indices \( S_1 \) where \( \sigma_{j}(t) \sigma_{j+1}(t) = 1 \). 
This case occurs if and only if the approximate function satisfies
\[
u_{j-\frac{1}{2}}(t) \leq u_{j+\frac{1}{2}}(t) \leq u_{j+\frac{3}{2}}(t) 
\quad \text{or} \quad 
u_{j-\frac{1}{2}}(t) \geq u_{j+\frac{1}{2}}(t) \geq u_{j+\frac{3}{2}}(t).
\]

\item The set of indices \( S_2 \) where \( \sigma_{j}(t) \sigma_{j+1}(t) = -1 \). 
This case occurs if and only if \( u_{j+\frac{1}{2}}(t) \) is 
a minimum or maximum value of the set given by
\(\{u_{j-\frac{1}{2}}(t), u_{j+\frac{1}{2}}(t), u_{j+\frac{3}{2}}(t)\}\). 
In this case, by using the definition of the slope limiter given in 
\eqref{slope}, and the definition of the min-mod function given in 
\eqref{minmod}, we can conclude that:
\[
s_{j+\frac{1}{2}}(t) = 0.
\]
\end{enumerate}
Therefore, in \eqref{diff-2}, the terms where \( \sigma_{j+1}(t) \cdot \sigma_{j}(t) = 1 \) 
vanish, and the terms where \( \sigma_{j+1}(t) \cdot \sigma_{j}(t) = -1 \) reduce to
\begin{align}
&\sum_{j}\frac{d}{dt}|\Delta u_{j+1}(t)|\label{diff-3}\\
&=\frac{2}{\Delta x}\sum_{j\in S_2}\Big\{A_{j+1}(t)|\Delta u_{j+1}(t)|+B_{j}(t)|\Delta u_{j}(t)|\Big\}.\notag
\end{align}
By item 2), we have that \( u^{-}_{j+1}(t) = u^{+}_{j}(t) = u_{j+\frac{1}{2}} \). 
Recalling the definitions of the terms \( A_{j+1}(t) \) and \( B_j(t) \) given in 
\eqref{A} and \eqref{B}, we have, respectively
\begin{align}
A_{j+1}(t)&=\frac{\partial }{\partial v}G(u^{-}_{j+1}(t),u^*)\Big(\frac{u^{+}_{j+1}(t)-u^{+}_{j}(t)}{u_{j+\frac{3}{2}}(t)-u_{j+\frac{1}{2}}(t)}\Big)\notag\\
&=\frac{\partial }{\partial v}G(u^{-}_{j+1}(t),u^*)\Big(1-\frac{\Delta x}{4}\frac{s_{j+\frac{3}{2}}}{\Delta u_{j+1}(t)}\Big)\label{part-1}
\end{align}
for some $u^*\in \big[\min\lbrace  u^{+}_{j+1}(t), u^{+}_{j}(t) \rbrace,\max\lbrace u^{+}_{j+1}(t), u^{+}_{j}(t) \rbrace \big]$ and
\begin{align}
B_{j}(t)&=\frac{\partial }{\partial u}G(u^{**},u^{+}_{j}(t))\Big(\frac{u^{-}_{j+1}(t)-u^{-}_{j}(t)}{\Delta u_{j}(t)}\Big)\notag\\
&=\frac{\partial }{\partial u}G(u^{**},u^{+}_{j}(t))\Big(1-\frac{\Delta x}{4}\frac{s_{j-\frac{1}{2}}}{\Delta u_{j}(t)}\Big),\label{part-2}
\end{align}
for some $u^{**}\in \big[\min\lbrace  u^{-}_{j+1}(t), u^{-}_{j}(t) \rbrace,\max\lbrace u^{-}_{j+1}(t), u^{-}_{j}(t) \rbrace \big]$.

Now, we can select the values in \eqref{part-1} and \eqref{part-2} 
according to the behavior of \( u_{j+\frac{1}{2}}(t) \) in the set  
\( \{u_{j-\frac{1}{2}}(t), u_{j+\frac{1}{2}}(t), u_{j+\frac{3}{2}}(t)\} \) as follows:
\begin{enumerate}[a).]
\item If \( u_{j+\frac{1}{2}}(t) \) is a maximum value of  
\( \{u_{j-\frac{1}{2}}(t), u_{j+\frac{1}{2}}(t), u_{j+\frac{3}{2}}(t)\} \), then we set
\[
A_{j+1}(t) = \frac{\partial }{\partial v}G(u^{-}_{j+1}(t), u^*)\Big(1 - \frac{\Delta x}{4}\frac{s_{j+\frac{3}{2}}}{\Delta u_{j+1}(t)}\Big) \leq 0,
\]
and
\[
B_{j}(t) = \frac{\partial }{\partial u}G(u^{**}, u^{+}_{j}(t))\Big(1 - \frac{\Delta x}{4}\frac{s_{j-\frac{1}{2}}}{\Delta u_{j}(t)}\Big) \leq 0.
\]

\item If \( u_{j+\frac{1}{2}}(t) \) is a minimum value of  
\( \{u_{j-\frac{1}{2}}(t), u_{j+\frac{1}{2}}(t), u_{j+\frac{3}{2}}(t)\} \), then we set
\[
A_{j+1}(t) = \frac{\partial }{\partial v}G(u^{-}_{j+1}(t), u^*)\Big(1 - \frac{\Delta x}{4}\frac{s_{j+\frac{3}{2}}}{\Delta u_{j+1}(t)}\Big) \leq 0,
\]
and
\[
B_{j}(t) = \frac{\partial }{\partial u}G(u^{**}, u^{+}_{j}(t))\Big(1 - \frac{\Delta x}{4}\frac{s_{j-\frac{1}{2}}}{\Delta u_{j}(t)}\Big) \leq 0.
\]
\end{enumerate}
Replacing the values \( A_{j+1}(t) \) and \( B_{j}(t) \) following the 
cases described in item a) and b) into \eqref{diff-3}, we obtain
\begin{equation}\label{TV-fin}
\sum_{j}\frac{d}{dt}|\Delta u_{j+1}(t)|\leq 0,
\end{equation} 
\textbf{Step 2.}
Now, notice that we can write 
\begin{align*}
& \frac{d}{dt}u_{\frac{1}{2}}(t)=-\frac{1}{\Delta x}A_{1}(t)\Delta u_{1}(t)+\frac{1}{\Delta x}B_{0}(t)\Delta u_{0}(t)\\
&+\frac{1}{\Delta x}(k_L-k_R)(1-u_{-\frac{1}{2}}(t))(u^+_0(t)+u^{-}_0(t)),
\end{align*}
and
\begin{align*}
& \frac{d}{dt}u_{-\frac{1}{2}}(t)=-\frac{1}{\Delta x}A_{-1}(t)\Delta u_{-1}(t)+\frac{1}{\Delta x}B_{-2}(t)\Delta u_{-2}(t)\\
&\frac{1}{\Delta x}(k_R-k_L)(1-u_{\frac{1}{2}}(t))(u^+_0(t)+u^{-}_0(t)).\\
\end{align*}

By using \eqref{TV-fin}, we have
\begin{align}
&\sum_{j}\frac{d}{dt}|\Delta u_{j+1}(t)|=\sum_{j>1}\frac{d}{dt}|\Delta u_{j+1}(t)|+\sum_{j<-1}\frac{d}{dt}|\Delta u_{j}(t)|+\sigma_{1}\frac{d}{dt}\Delta u_{1}(t)\notag\\
&+\sigma_{0}\frac{d}{dt}\Delta u_{0}(t)+\sigma_{-1}\frac{d}{dt}\Delta u_{-1}(t)\notag\\
&\leq (\sigma_0-\sigma_1)\frac{1}{\Delta x}(k_L-k_R)(1-u_{-\frac{1}{2}}(t))(u^+_0(t)+u^{-}_0(t))\notag\\
&+(\sigma_{-1}-\sigma_0)\frac{1}{\Delta x}(k_R-k_L)(1-u_{\frac{1}{2}}(t))(u^+_0(t)+u^{-}_0(t))\notag\\
&\leq \frac{1}{\Delta x}\sigma_0(t)(k_L-k_R)(u^+_0(t)+u^{-}_0(t))\times\notag\\
&\Big(\big(1-\sigma_0(t)\sigma_1(t)\big)(1-u_{-\frac{1}{2}}(t))+\big(1-\sigma_0(t)\sigma_{-1}(t)\big)(1-u_{\frac{1}{2}}(t))\Big)\label{TV-fin-2}
\end{align}
Letting 
\begin{align*}
B(t)&=\sigma_0(t)(k_L-k_R)(u^+_0(t)+u^{-}_0(t))\\ &\times\Big(\big(1-\sigma_0(t)\sigma_1(t)\big)(1-u_{-\frac{1}{2}}(t))+\big(1-\sigma_0(t)\sigma_{-1}(t)\big)(1-u_{\frac{1}{2}}(t)\Big),
\end{align*}
this function  is bounded by Proposition \ref{p-maximum}. 
Then, applying the definition of total variation given 
in \eqref{TV-norm} and  using  \eqref{TV-fin-2}, we obtain 
{
\begin{equation*}
\sum_{j}|\Delta u_{j+1}(t)|=TV(u_0)+\frac{1}{\Delta x}\int_{0}^{t}B(s)ds.
\end{equation*}
}
\end{proof}
\begin{rmk}
Notice that inequality \eqref{TV-b} from Proposition \ref{TV-nounif} 
does not provide a uniform bound on the total variation of the 
sequence \( \{u_{j+\frac{1}{2}}(t)\}_{j\in \mathbb{Z}} \). This fact 
is consistent with the behavior of solutions in the context of 
scalar conservation laws with discontinuous flux. Furthermore, 
this result aligns with the BV-regularity issues of AB-entropy solutions in \cite{Adimurthi1}, Theorem 2.13.
On the other hand, we have a uniform bound in \( BV_{\text{loc}}(\mathbb{R} \setminus \{ 0 \}) \), as expected considering that we seek for classic week entropy solutions outside the discontinuity.

\end{rmk}
\subsection{Temporal continuity}

In order to analyze the temporal continuity, let us recall the 
ODE \eqref{Semi-pm}.
\begin{prop}\label{cont-modulus}
Let \( u^{\Delta x} \) be a solution to the initial value problem 
\eqref{Cauchy-1}-\eqref{Cauchy-Data}. Suppose that \( u_0 \in BV(\mathbb{R}) \). 
Then, for any \( \tau > 0 \), there exist two positive constants \( K \) and \( L \), such that 
\begin{equation}\label{TC-01}
\|u^{\Delta x}(t+\tau,\cdot) - u^{\Delta x}(t,\cdot)\|_{L^1(\mathbb{R})} \leq  \tau K \Big(TV(u_{0}) + L\Big),
\end{equation}
for all \( t \in [0, T-\tau] \).
\end{prop}

\begin{proof}
Since $u^{\Delta}$ is a solution to the ODE \eqref{Cauchy-1},  we have that
\begin{align}
&\frac{d}{dt}u_{j+\frac{1}{2}}(t)=-\frac{1}{\Delta x}\Big\{\Big(\mathcal{G}_{j+1}(u^{-}_{j+1}(t),u^{+}_{j+1}(t))-\mathcal{G}_{j}(u^{-}_{j}(t),u^{+}_{j+1}(t))\Big).\label{TC-1}
\end{align}
In \eqref{TC-01}, recalling the definition of \( \mathcal{G}_{j} \) given in \eqref{G-func}, we have
\begin{align}
&\mathcal{G}_{j+1}(u,v)-\mathcal{G}_{j}(r,s)=\frac{1}{4}\big(b_{j+1}(t)(u-v)-b_{j}(t)(r-s))\notag\\
&+\frac{1}{4}\Big(\big(f_{j+\frac{1}{2}}+f_{j+\frac{3}{2}}\big)(u+v)-\big(f_{j-\frac{1}{2}}+f_{j+\frac{1}{2}}\big)(r+s)\Big).\label{TC-2}
\end{align}
From \eqref{TC-2}, the behavior of \( \mathcal{G}_{j+1}(u,v) - \mathcal{G}_{j}(r,s) \) 
depends on the value of the index \( j \). Specifically, we have:

\begin{enumerate}[1).]
\item For \( j < -1 \), we have:
\begin{align}
\mathcal{G}_{j+1}(u,v) - \mathcal{G}_{j}(r,s) &= 
\frac{1}{4} \big(b_{j+1}(t)(u-v) - b_{j}(t)(r-s)\big) \label{TCdif-1} \\
&\quad + \frac{1}{4} k_L \big((2-(u+v))(u+v) - (2-(r+s))(r+s)\big). \notag
\end{align}

\item For \( j = -1 \), we have:
\begin{align}
\mathcal{G}_{0}(u,v) - \mathcal{G}_{-1}(r,s) &=
\frac{1}{4} \big(b_{0}(t)(u-v) - b_{-1}(t)(r-s)\big) \label{TCdif-2} \\
&\quad + k_L \big((1-u)(u+v) - (1-r)(r+s)\big) \notag \\
&\quad + k_R \big((1-v)(u+v) - (1-s)(r+s)\big) \notag \\
&\quad + \frac{1}{4}(k_L - k_R)(1-v)(u+v). \notag
\end{align}

\item For \( j = 0 \), we have:
\begin{align}
\mathcal{G}_{1}(u,v) - \mathcal{G}_{0}(r,s) &=
\frac{1}{4} \big(b_{1}(t)(u-v) - b_{0}(t)(r-s)\big) \label{TCdif-3} \\
&\quad + k_L \big((1-u)(u+v) - (1-r)(r+s)\big) \notag \\
&\quad + k_R \big((1-v)(u+v) - (1-s)(r+s)\big) \notag \\
&\quad + \frac{1}{4}(k_R - k_L)(1-v)(u+v). \notag
\end{align}

\item For \( j > 0 \), we have:
\begin{align}
\mathcal{G}_{j+1}(u,v) - \mathcal{G}_{j}(r,s) &= 
\frac{1}{4} \big(b_{j+1}(t)(u-v) - b_{j}(t)(r-s)\big) \label{TCdif-4} \\
&\quad + \frac{1}{4} k_R \big((2-(u+v))(u+v) - (2-(r+s))(r+s)\big). \notag
\end{align}
\end{enumerate}

Let
\[
g(u,v) = (2 - (u+v))(u+v), \quad g_{-}(u,v) = (1-u)(u+v), \quad g_{+}(u,v) = (1-v)(u+v),
\]
these functions are locally Lipschitz. Therefore, the terms on the 
left-hand side of \eqref{TCdif-1}, \eqref{TCdif-2}, \eqref{TCdif-3}, and \eqref{TCdif-4} 
are composed of locally Lipschitz functions as well. Set
\[
\text{L}_g = \max\{\text{Lip}(g), \text{Lip}(g_-), \text{Lip}(g_+)\},
\]
then, replacing in \eqref{TCdif-1}, \eqref{TCdif-2}, \eqref{TCdif-3}, 
and \eqref{TCdif-4}, we obtain:
\begin{enumerate}[1).]
\item For \( j < -1 \):
\begin{align}
&|u_{j+\frac{1}{2}}(t+\tau) - u_{j+\frac{1}{2}}(t)| \label{Tsum-1} \\
&\leq \int_{t}^{t+\tau} \Big(b_{j+1}(s)|u^{+}_{j+1}(s) - u^{-}_{j+1}(s)| 
+ b_{j}(s)|u^{+}_{j}(s) - u^{-}_{j}(s)|\Big) ds \notag \\
&\quad + k_L \text{L}_{g} \int_{t}^{t+\tau} \Big(|u^{+}_{j+1}(s) - u^{-}_{j+1}(s)| 
+ |u^{+}_{j}(s) - u^{-}_{j}(s)|\Big) ds. \notag
\end{align}
\item For $j=-1$
\begin{align}
&|u_{-\frac{1}{2}}(t+\tau)-u_{-\frac{1}{2}}(t)|\notag\\
&\leq \int_{t}^{t+\tau}\Big(b_{0}(s)|u^{+}_{0}(s)-u^{-}_{0}(s)|+b_{-1}(s)|u^{+}_{-1}(s)-u^{-}_{-1}(s)|\Big)ds\notag\\
&+k_L \text{L}_{g}\int_{t}^{t+\tau}\Big(|u^{+}_{0}(s)-u^{-}_{0}(s)|+|u^{+}_{-1}(s)-u^{-}_{-1}(s)|\Big)ds\notag\\
&+k_R \text{L}_{g}\int_{t}^{t+\tau}\Big(|u^{+}_{0}(s)-u^{-}_{0}(s)|+|u^{+}_{-1}(s)-u^{-}_{-1}(s)|\Big)ds.\label{Tsum-2}
\end{align}
\item For $j=0$
\begin{align}
&|u_{\frac{1}{2}}(t+\tau)-u_{\frac{1}{2}}(t)|\notag\\
&\leq \int_{t}^{t+\tau}\Big(b_{1}(s)|u^{+}_{1}(s)-u^{-}_{1}(s)|+b_{0}(s)|u^{+}_{0}(s)-u^{-}_{0}(s)|\Big)ds\notag\\
&+k_L \text{L}_{g} \int_{t}^{t+\tau}\Big(|u^{+}_{1}(s)-u^{-}_{1}(s)|+|u^{+}_{0}(s)-u^{-}_{0}(s)|\Big)ds\notag\\
&+k_R \text{L}_{g} \int_{t}^{t+\tau}\Big(|u^{+}_{1}(s)-u^{-}_{1}(s)|+\text{L}^2_{G}|u^{+}_{0}(s)-u^{-}_{0}(s)|\Big)ds.\label{Tsum-3}
\end{align}
\item For $j>0$
\begin{align}
&|u_{j+\frac{1}{2}}(t+\tau)-u_{j+\frac{1}{2}}(t)|\notag\\
&\leq \int_{t}^{t+\tau}\Big(b_{j+1}(s)|u^{+}_{j+1}(s)-u^{-}_{j+1}(s)|+b_{j}(s)|u^{+}_{j}(s)-u^{-}_{j}(s)|\Big)ds\notag\\
&+k_R \text{L}_{g} \int_{t}^{t+\tau}\Big(|u^{+}_{j+1}(s)-u^{-}_{j+1}(s)|+|u^{+}_{j+1}(s)-u^{-}_{j+1}(s)|\Big)ds.\label{Tsum-4}
\end{align}
\end{enumerate}
Finally we set
\[
K=\text{L}_g\cdot\max\{k_L, k_R\}.
\]
Replacing $K$ into \eqref{Tsum-1}, \eqref{Tsum-2}, \eqref{Tsum-3} and \eqref{Tsum-4}, summing over $j$ and taking  into account \eqref{Tsum-1}, \eqref{Tsum-2}, \eqref{Tsum-3} and \eqref{Tsum-4}, we obtain
\begin{align}
&\sum_{j}|u_{j+\frac{1}{2}}(t+\tau)-u_{j+\frac{1}{2}}(t)|\Delta x\label{TC-1-new}\\
&\leq \int_{t}^{t+\tau}\sum_{j}\Big(b_{j+1}(s)|u^{+}_{j+1}(s)-u^{-}_{j+1}(s)|+b_{j}(s)|u^{+}_{j}(s)-u^{-}_{j}(s)|\Big)\Delta x ds\notag\\
&+K \int_{t}^{t+\tau}\sum_{j}\Big(|u^{+}_{j+1}(s)-u^{-}_{j+1}(s)|+|u^{+}_{j+1}(s)-u^{-}_{j+1}(s)|\Big)\Delta x ds.\notag
\end{align}
Now, by using the definition of $u^{\pm}_{j}(t)$ given in \eqref{upm}, we get
\begin{equation}\label{Mod-cont}
|u^{+}_{j+1}(t+\tau)-u^{-}_{j+1}(t)|\leq \frac{3}{4}|u_{j+\frac{3}{2}}(t)-u_{j+\frac{1}{2}}(t)|
\end{equation}
By Proposition \ref{p-maximum}, the sequence of functions 
\( b_{j}(t) \) is uniformly bounded. If we set \( b = \sup_{j}\{b_j(t)\} \), 
replacing \eqref{Mod-cont} into \eqref{TC-1-new} and using \ref{TV-nounif}, we obtain
\[
\sum_{j} |u_{j+\frac{1}{2}}(t+\tau) - u_{j+\frac{1}{2}}(t)| \Delta x 
\leq \tau (K + b) \Big(TV(u_{0}) + \int_{0}^{t} |B(s)| ds\Big).
\]
Finally, we set \( L = \sup_{t \in [0,T]} \int_{0}^{t} |B(s)| ds \) to obtain \eqref{TC-01}.
\end{proof}

\subsection{Convergence to  entropic solutions}\label{cwesol}
In this section, we aim to show that the approximations produced by the semidiscrete Lagrangian-Eulerian scheme converge to the entropic solution for \eqref{1-1}-\eqref{1-1-D} in the sense of Definition \eqref{3-loc-weak}.
Recalling the flux function defined in \eqref{num_flux}, we denote
\begin{align}
&\mathcal{F}_j(t)=\mathcal{F}(u^-_j(t),u^+_j(t))=\label{ent24-1}\\
&\frac{1}{4}\Big[b_{j}(t)\Big(u^{-}_{j}(t)-u^+_{j}(t)\Big)+\big(f_{j-\frac{1}{2}}(t)+f_{j+\frac{1}{2}}(t)\big)\Big(u^-_{j}(t)+u^+_{j}(t)\Big)\Big],\notag
\end{align}
and rewrite the semidiscrete scheme in the form 
\begin{equation}\label{ent24-2}
u^{\prime}_{j+\frac{1}{2}}(t)=-\frac{1}{\Delta x}\Big(\mathcal{F}_{j+1}(t)-\mathcal{F}_{j}(t)\Big),
\end{equation}
then we can state the following result.

\begin{thm}
Let \( u^{\Delta} \) be defined in \eqref{Fin-app} and obtained by the semi-discrete Lagrangian-Eulerian scheme \eqref{Semi}. Then, \( u^{\Delta} \to u \) in \( L^1_{\text{loc}}([0,T] \times (\mathbb{R}\setminus\lbrace 0\rbrace)) \), where \( u \) is the unique entropy solution of the Cauchy problem \eqref{1-1}, \eqref{1-1-D}.
\end{thm}

\begin{proof}By Propositions \ref{p-maximum}, 
$$  \Vert u^{\Delta} \Vert_{L^\infty}\leq 1 \quad \hbox{for all } t\in[0,T], \hbox{for all }n\in\mathbb{N}. $$
Hence $u^{\Delta}$ is $\emph{weakly}^*$ relatively compact in $L^\infty([0,T]\times\mathbb{R})$.
By \ref{L1-est} and \ref{TV-nounif}, it follows that $u^{\Delta}\in L^1(\mathbb{R})\cap BV_{loc}(\mathbb{R}\setminus\lbrace 0 \rbrace)$ and the total variation is uniformly bounded in the compact sets of $\mathbb{R}\setminus\lbrace 0 \rbrace$. We also have, by \ref{cont-modulus}, that $u^{\Delta}$ is locally Lipschitz in time with respect to the $L^1$ norm.  Therefore, by Helly's compactness theorem, there exists a subsequence, still labeled \( u^{\Delta} \), and a function \( u \in L^1_{\text{loc}}([0,T] \times(\mathbb{R}\setminus\lbrace 0\rbrace)) \) such that \( u^{\Delta} \to u \) in \( L^1_{\text{loc}}([0,T] \times (\mathbb{R}\setminus\lbrace 0\rbrace)) \). \\
Now, we show that \( u \) satisfies the entropy inequality \eqref{3-loc-weak}. Let 
\[
H_{j}(t) = \mathcal{F}_j(u^-_j(t), u^+_j(t)) - \mathcal{F}_j(c, c),
\]
where \( c \in [0,1] \).

We rewrite \eqref{ent24-2} in the following form
\begin{equation}\label{Pro-ent-1}
u^{\prime}_{j+\frac{1}{2}}(t)=
\begin{cases}
-\frac{1}{\Delta x}\Big(H_{j+1}(t)-H_{j}(t) \Big) \hbox{ if }j<-1 \hbox{ or }j\geq 1\\
-\frac{1}{\Delta x}\Big(H_{0}(t)-H_{-1}(t) \Big)+\frac{1}{2\Delta x} (k_L-k_R)f(c)\hbox{ if } j=-1\\
-\frac{1}{\Delta x}\Big(H_{1}(t)-H_{0}(t) \Big)+\frac{1}{2\Delta x} (k_L-k_R)f(c)\hbox{ if } j=0

\end{cases}
\end{equation}
Let \( \eta \) be a smooth convex function. We can assume without loss of generality that \( |\eta'| \leq 1 \). Multiplying by \( \eta^{\prime} \) in \eqref{Pro-ent-1} when \( j \neq 0, -1 \), we obtain:
\[
\frac{d}{dt} \eta(u_{j+\frac{1}{2}}(t)) = -\frac{1}{\Delta x} \eta^{\prime}(u_{j+\frac{1}{2}}(t)) \Big(H_{j+1}(t) - H_{j}(t)\Big).
\]
By using the notation
\[
D_+\Big(\eta^{\prime}(u_{j-\frac{1}{2}}(t))H_j(t)\Big)=\eta^{\prime}(u_{j+\frac{1}{2}}(t))H_{j+1}(t)+\eta^{\prime}(u_{j-\frac{1}{2}}(t)) H_j(t),
\]
we obtain the relation
\begin{align}
&\eta^{\prime}(u_{j+\frac{1}{2}}(t))=\notag\\
&-\frac{1}{\Delta x}\Big(D_+\Big(\eta^{\prime}(u_{j-\frac{1}{2}}(t))H_j(t)\Big)-D_+\eta^{\prime}(u_{j-\frac{1}{2}}(t)) H_j(t)\Big).\label{Proof-ent-2}
\end{align}
Given $\phi \in C^1_c([0,\infty)\times \mathbb{R})$ non negative, we  define $\phi_{j+\frac{1}{2}}(t)=\phi(t,x_{j+\frac{1}{2}})$. 
Multiplying \eqref{Proof-ent-2} by $\phi_{j+\frac{1}{2}}\Delta x$, we obtain
\begin{align}
&\frac{d}{dt}\eta(u_{j+\frac{1}{2}}(t))\phi_{j+\frac{1}{2}}(t)d\Delta x=\notag\\
&-\frac{1}{\Delta x}D_+\Big(\eta^{\prime}(u_{j-\frac{1}{2}}(t))H_j(t)\Big)\phi_{j+\frac{1}{2}}(t)\Delta x\label{Int-24-1}\\
&+\frac{1}{\Delta x}D_+\eta^{\prime}(u_{j-\frac{1}{2}}(t)) H_j(t)\phi_{j+\frac{1}{2}}(t)\Delta x.\notag
\end{align}
Replacing the  identity
\begin{align*}
&D_+\Big(\eta^{\prime}(u_{j-\frac{1}{2}}(t))H_j(t)\phi_{j-\frac{1}{2}}(t)\Big)\\
&=D_+\Big(\eta^{\prime}(u_{j-\frac{1}{2}}(t))(t)H_j(t)\Big)\phi_{j+\frac{1}{2}}(t)+\eta^{\prime}(u_{j-\frac{1}{2}}(t))H_j(t)D_+\Big(\phi_{j-\frac{1}{2}}(t)\Big),
\end{align*}
into \eqref{Int-24-1}, we obtain
\begin{align}
&\frac{d}{dt}\eta(u_{j+\frac{1}{2}}(t))\phi_{j+\frac{1}{2}}(t)d\Delta x=\eta^{\prime}(u_{j-\frac{1}{2}}(t))H_j(t)\frac{1}{\Delta x}D_+\Big(\phi_{j-\frac{1}{2}}(t)\Big)\Delta x\notag\\
&-\frac{1}{\Delta x}D_+\Big(\eta^{\prime}(u_{j-\frac{1}{2}}(t))H_j(t)\phi_{j-\frac{1}{2}}(t)\Big)\Delta x\label{Int-24-3}\\
&+\frac{1}{\Delta x}D_+\eta^{\prime}(u_{j-\frac{1}{2}}(t)) H_j(t)\phi_{j+\frac{1}{2}}(t)\Delta x.\notag
\end{align}
Now performing the same computation in the case of $j=0,-1$, we get
\begin{align}
&\frac{d}{dt}\eta(u_{j+\frac{1}{2}}(t))\phi_{j+\frac{1}{2}}(t)d\Delta x=\eta^{\prime}(u_{j-\frac{1}{2}}(t))H_j(t)\frac{1}{\Delta x}D_+\Big(\phi_{j-\frac{1}{2}}(t)\Big)\Delta x\notag\\
&-\frac{1}{\Delta x}D_+\Big(\eta^{\prime}(u_{j-\frac{1}{2}}(t))H_j(t)\phi_{j-\frac{1}{2}}(t)\Big)\Delta x\label{Int-24-3}\\
&+\frac{1}{\Delta x}D_+\eta^{\prime}(u_{j-\frac{1}{2}}(t)) H_j(t)\phi_{j+\frac{1}{2}}(t)\Delta x +\frac{1}{2 }\eta^{\prime}(u_{j+\frac{1}{2}}(t)) (k_L-k_R)f(c)\phi_{j+\frac{1}{2}}\notag
\end{align}

We sum over $j$ on both the left and right-hand sides of \eqref{Int-24-3}, and since $\phi$ is a function of compact support, this sum is finite. Taking into account the terms for $j = -1$ and $j = 0$, we obtain
\begin{align}
&\sum_{j}\frac{d}{dt}\eta(u_{j+\frac{1}{2}}(t))\phi_{j+\frac{1}{2}}(t)\Delta x=\sum_{j}\eta^{\prime}(u_{j-\frac{1}{2}}(t))H_j(t)\frac{1}{\Delta x}D_+\phi_{j+\frac{1}{2}}(t)\Delta x\notag\\
&+\sum_{j}\frac{1}{\Delta x}D_+\eta^{\prime}(u_{j-\frac{1}{2}}(t)) H_j(t)\phi_{j+\frac{1}{2}}(t)\Delta x\label{Int-24-4}\\
&+(k_L-k_R)\frac{\Big(\eta^{\prime}(u_{+\frac{1}{2}}(t)) \phi_{+\frac{1}{2}}(t)+\eta^{\prime}(u_{-\frac{1}{2}}(t)) \phi_{-\frac{1}{2}}(t)\Big)}{2}f(c)\label{ext-term}
\end{align}
By the assumption $|\eta'|\leq 1$,  we have the inequality
\[
\eqref{ext-term}\leq|k_R-k_L|\frac{\phi_{+\frac{1}{2}}(t)+\phi_{-\frac{1}{2}}(t)}{2}f(c).
\]
The term in \eqref{Int-24-4} tends to 0 as $\Delta x \to 0$. Taking the limit as $\Delta x \to 0$ and integrating over the temporal variable, we obtain
\begin{align*}
&-\int_{R}\eta(u(0,x))\partial_t \phi(0,x)-\int_{0}^{\infty}\eta(u(t,x))\partial_t \phi(t,x)dxdt\\
&\leq \int_{0}^{\infty}\int_{R} k(x)\eta^{\prime}(u(t,x))\big(f(u(t,x)-f(c)))\partial_x\phi(t,x) +|k_R-k_L|f(c)\int_{0}^{\infty}\phi(t,0)dt.\notag
\end{align*}
Finally, we consider a sequence $\eta_r$ such that $\eta_r \to |\cdot|$ and $\eta_r^{\prime}(r) \to \text{sgn}(\cdot)$. Replacing $\eta$ by $\eta_r(u - c)$ and taking the limit as $r \to 0$, we obtain
\begin{align*}
&\int_{0}^{\infty}\int_{\mathbb{R}}\Big(|u(t,x)-c|\partial_t \phi(t,x)+k(x)\Phi(u(t,x),c)\partial_x\phi(t,x)\Big)dt\,dx \notag\\
&+\int_{\mathbb{R}}{|u_0(x)-c|}\phi(0,x)dx+|k_{R}-k_{L}|f(c)\int_{0}^{\infty}\phi(t,0)dt\geq 0
\end{align*}
\end{proof}

\section{Numerical examples}\label{numexa}
\subsection{Examples of local flux with a discontinuous flux}
In this last section, we apply our scheme to several examples already discussed in the existing literature (\cite[Section 5]{Vovelle}, \cite[Section 5.2]{Kenneth2}, and \cite[Section 7]{Shen-1}). For each of these examples, we will compare the approximate solutions $u_{\text{app}}$ with the exact ones $u_{\text{exact}}$. We will compute their $L^1$-distance using the $L^1$-norm defined in \eqref{L1-norm}, the relative error defined by $\displaystyle \frac{\|u_{\text{exact}} - u_{\text{app}}\|_{L^1}}{\|u_{\text{exact}}\|_{L^1}}$, and the total variation of $u_{\text{app}}$, defined in \eqref{TV-norm}, taking cells of order $2^k + 1$.
\\

\begin{exmp}\label{Example-5.1}
In this example we reproduce the results of Seguin-Vovelle \cite[Section 5]{Vovelle}. 
In \eqref{k-f}, we consider the function
\begin{equation}\label{k-vovelle}
k(x)=H(x)+2(1-H(x)).
\end{equation}
and analyse the solutions of the Cauchy problem \eqref{1-1} with initial datum
\begin{equation}\label{Dat-ex1}
u_0(x)=\begin{cases}
0.5, \quad \text{if}\ x<0,\\
0.3, \quad \text{if}\ x>0.\\
\end{cases}
\end{equation}
In this case, $k_L > k_R$, the left and right states in the Riemann problem are given by $u_L = 0.5$ and $u_R = 0.3$. The entropic solution presents a shock wave connecting $u_L$ to $u^-$, where $u^- = \frac{1}{2} \left( 1 + \frac{1}{\sqrt{2}} \right)$, and a rarefaction wave connecting $u^+$ to $u_R$, with $u^+ = \frac{1}{2}$, i.e.
\begin{equation}\label{Ex1-sol}
u_{exact}(t,x)=\begin{cases}
\frac{1}{2},\quad \text{if} \ x< -\frac{1}{\sqrt{2}}t,\\
\frac{1}{2}\Big(1+\frac{1}{\sqrt{2}}\Big) ,\quad \text{if} \ -\frac{1}{\sqrt{2}}t<x <0,\\
\frac{1}{2}\Big(1-\frac{x}{t}\Big), \quad \text{if}\ 0\leq x\leq \frac{2}{5}t,\\
\frac{3}{10}, \quad \text{if}\ x> \frac{2}{5}t.
\end{cases}
\end{equation}
In Figure \ref{fig:1}, we illustrate our results by comparing the exact solution \eqref{Ex1-sol} with the approximate solutions at times $t = 0.5$ and $t = 1$.
\begin{figure}[h!]
\begin{center}
\begin{tabular}{c c}
\includegraphics[scale=0.35]{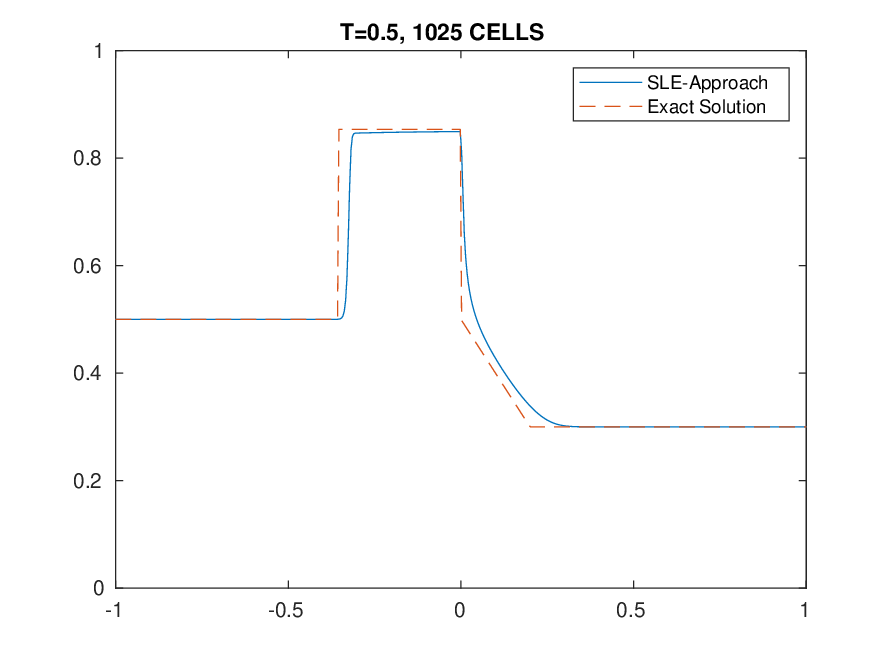} &
\includegraphics[scale=0.35]{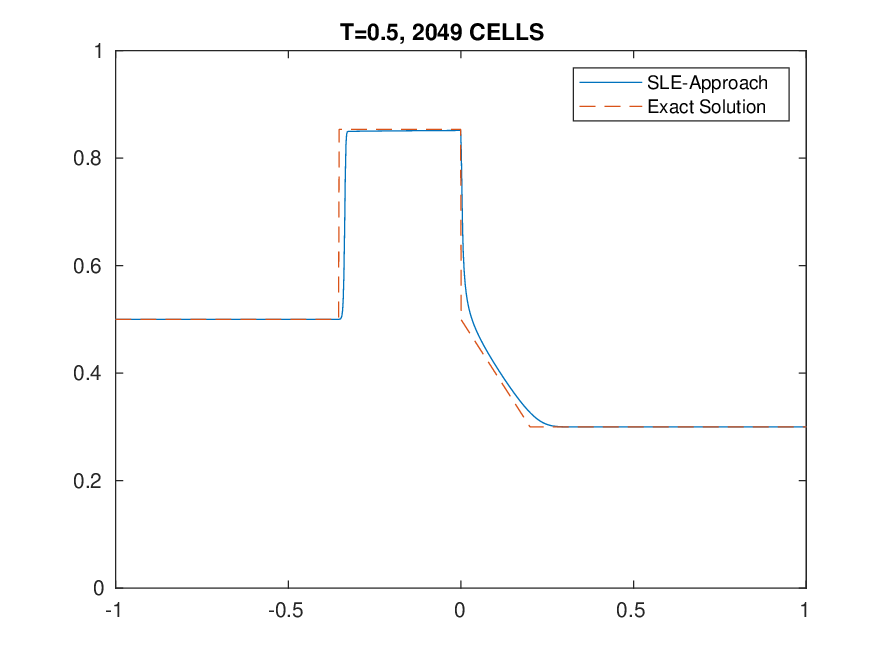} \\
($a$) & ($b$) \\
\includegraphics[scale=0.35]{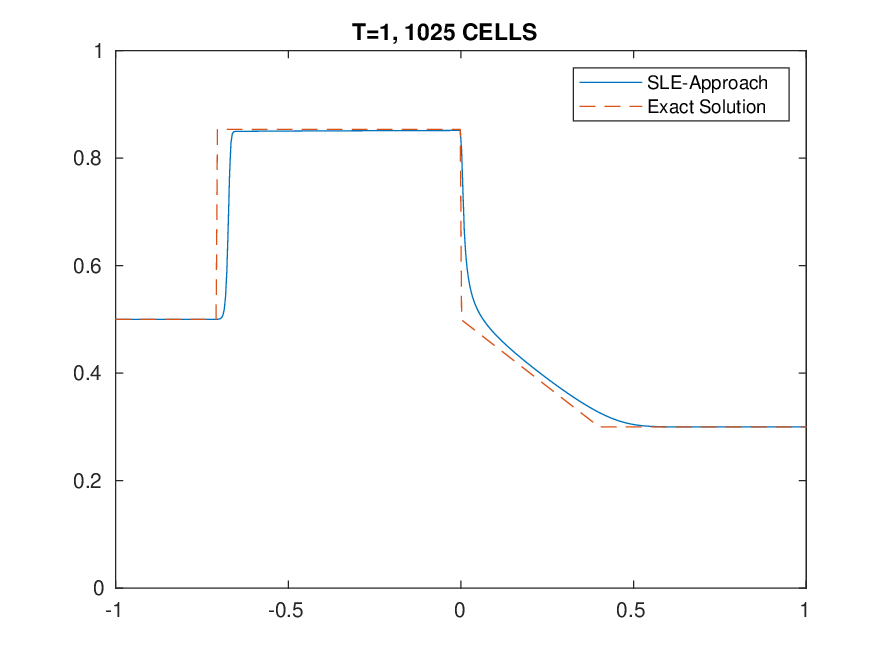} &
\includegraphics[scale=0.35]{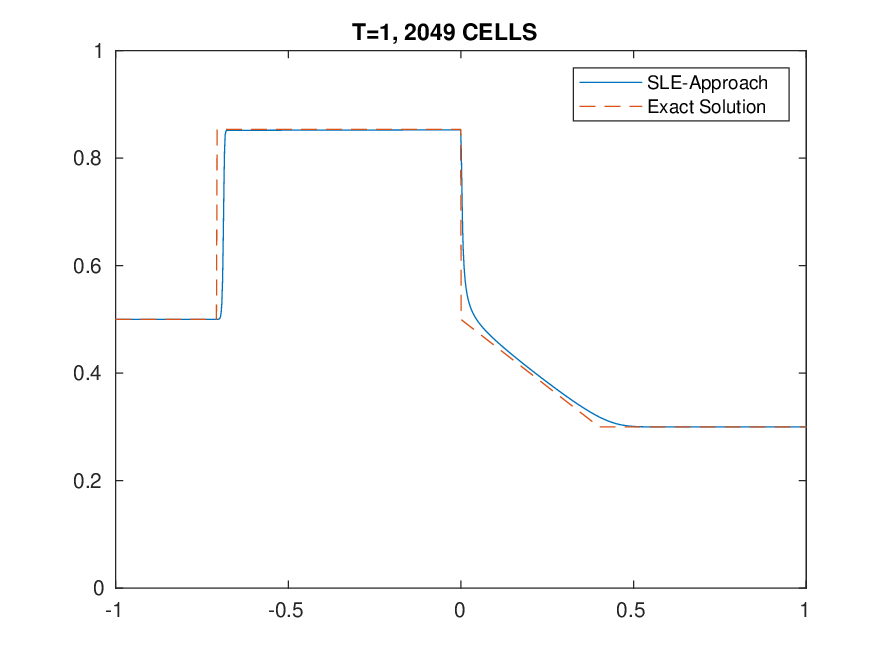} \\
($c$) & ($d$)
\end{tabular}
\end{center}
\caption{Illustration  comparing the exact solution of the 
Cauchy problem \eqref{1-1} with datum \eqref{Dat-ex1}. The exact solution 
is given by  \eqref{Ex1-sol} with the approximate solutions. 
a) and b) at time $t=0.5$ considering 1025 and 2049 cells. 
c) and d) results al time $t=1$ with 1025 and 2049 cells.
}
\label{fig:1}
\end{figure}
In Table \ref{Table:1}, we provide the $L^1$ error, the relative error, and the total variation  of the approximate solutions.
\begin{table}[h!]
    \centering
\begin{tabular}{ |c|c|c|c| } 
 \hline
 cells & $L^1$-error & Error Relative & $TV$ approximate solution \\ \hline
 257 & 0.0656 &0.0377 &0.8777 \\ 
 513 & 0.0383 &0.0220 &0.8915 \\
 1025 & 0.0240 &0.0137 &0.8990 \\ 
 2049 & 0.0137 &0.0078 &0.9030 \\
 4097 & 0.0078 &0.0045 &0.9050 \\
 \hline
\end{tabular}

    \caption{$L^1$-error, $L^1$-relative error and and total variation of approximate  solutions for Example \ref{Example-5.1}}
    \label{Table:1}
\end{table}

\end{exmp}
\newpage
\begin{exmp}\label{Example-5.5}
Following \cite{Vovelle}, in this example, we consider the same function $k$ given in \eqref{k-vovelle} and the initial datum
\begin{equation}\label{Dataex-2}
u_0(x)=\begin{cases}
0.95, \quad \text{if}\ x<0,\\
0.8, \quad \text{if}\ x>0.
\end{cases}
\end{equation}
Since in this case we have $u_R > 0.5$, from \cite[Appendix A.2]{Vovelle}, we obtain $u^- = \frac{1}{2} \left( 1 + \frac{\sqrt{17}}{5} \right)$ and $u^+ = u_R$. The solution consists of a rarefaction wave connecting $(u_L, u^-)$ and a shock wave connecting $(u^+, u_R)$.
\begin{equation}\label{Ex2-sol}
u_{exact}(t,x)=\begin{cases}
\frac{19}{20},\quad \text{if}\ x<-\frac{9}{5}t,\\
\frac{1}{4}\Big(2-\frac{x}{t}\Big),\quad \text{if}\ -\frac{9}{5}t<x<-2\frac{\sqrt{17}}{5}t,\\
\frac{1}{2}\Big(1+\frac{\sqrt{17}}{5}\Big),\quad \text{if}\  -2\frac{\sqrt{17}}{5}t<x<0,\\
\frac{4}{5}, \quad \text{if} \ x>0. 
\end{cases}
\end{equation}

\begin{figure}[ht!]
\begin{center}
\begin{tabular}{c c}
\includegraphics[scale=0.35]{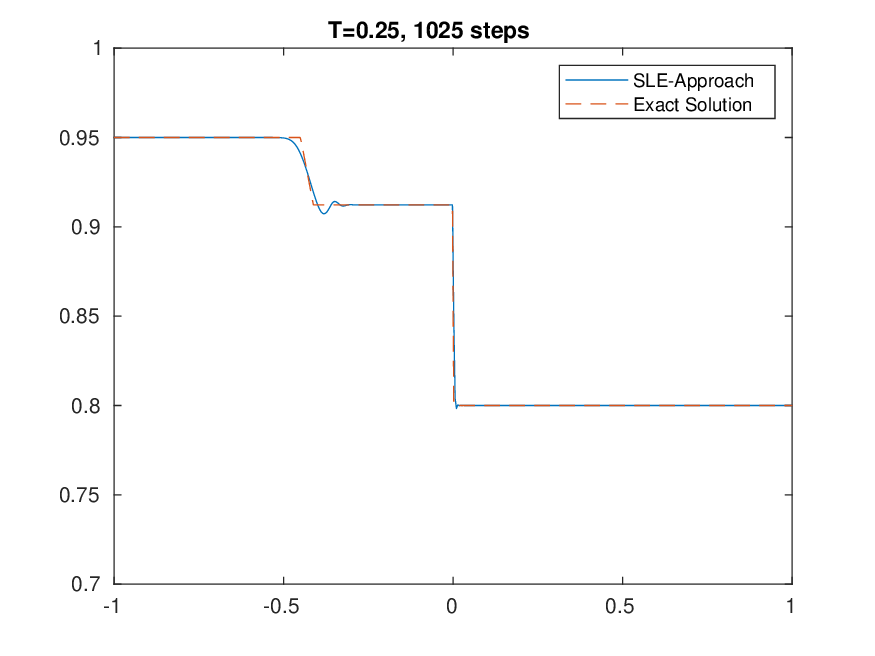} &
\includegraphics[scale=0.35]{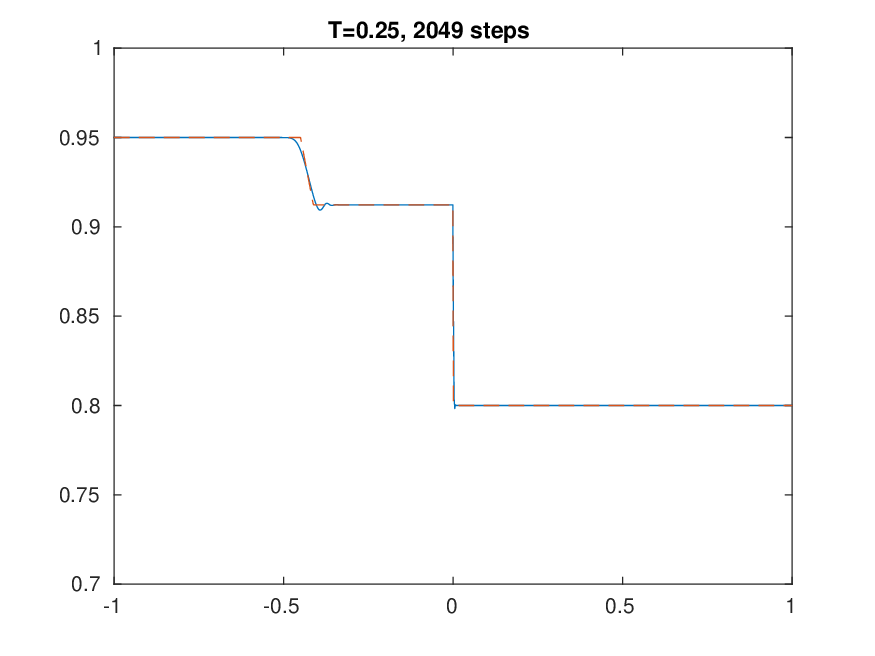} \\
($a$) & ($b$) \\
\includegraphics[scale=0.35]{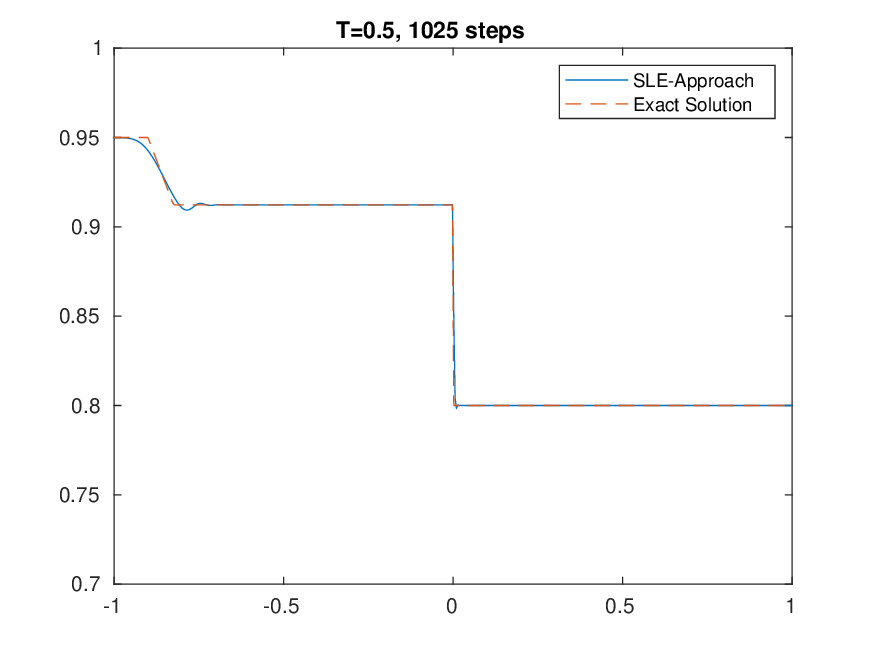} &
\includegraphics[scale=0.35]{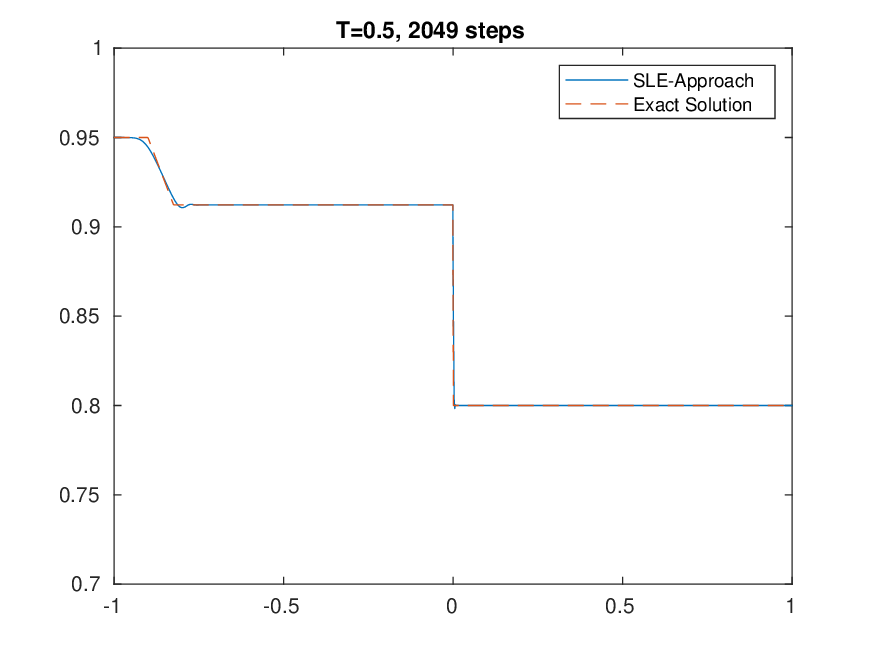} \\
($c$) & ($d$)
\end{tabular}
\end{center}
\caption{Illustration  comparing the exact solution of the 
Cauchy problem \eqref{1-1} with datum \eqref{Dataex-2}. The exact 
solution is given by \eqref{Ex2-sol}
}
\label{fig:2}
\end{figure}
In Table \ref{Table:2} we provide  again the $L^1$ error, the $L^1$-relative and the total variation  of the approximate solution.
\begin{table}[h!]
    \centering
\begin{tabular}{ |c|c|c|c| } 
 \hline
 cells & $L^1$-error & $L^1$-relative error &$TV$ approximate solution \\ \hline
 257 & 0.0029 & 0.0008259 &0.1778 \\ 
 513 & 0.0015 & 0.0004213&0.1697 \\
 1025 & 0.0008 &0.0002170 &0.1627 \\ 
 2049 & 0.0004 &0.0001142 &0.1582 \\
 4097 & 0.0002 &0.0000620 &0.1560 \\
 \hline
\end{tabular}

    \caption{$L^1$-error, $L^1$-relative error and total variation of approximate  solutions for Example \ref{Example-5.5}}
    \label{Table:2}
\end{table}

\end{exmp}


\begin{exmp}\label{Example-5.8}
As the last example, we consider a case in which $k_L < k_R$. Following \cite[Section 5]{Kenneth2}, we consider
\begin{equation}\label{k-keneth}
k(x)=H(x)+\frac{1}{2}(1-H(x)),
\end{equation}
and as initial datum
\begin{equation}\label{Exm3-sol}
u(t,x)=\begin{cases}
0.8, \quad \text{if} \ x<0,\\
0.1, \quad \text{if} \ x<0.
\end{cases}
\end{equation}
The exact solution presents a shock wave connecting $(u_L, u^-)$ with $u^- = \frac{1}{2} \left( 1 + \frac{1}{\sqrt{2}} \right)$ and a rarefaction wave connecting $(u^+, u_R)$ with $u^+ = \frac{1}{2} \left( 1 - 2 \frac{x}{t} \right)$, i.e.
\begin{equation}\label{Sol-3}
u_{exact}(t,x)=\begin{cases}
\frac{4}{5},\quad \text{if}\  x<\frac{1}{2}\Big(\frac{1}{2\sqrt{2}}+\frac{4}{5}\Big)t,\\
\frac{1}{2}\Big(1+\frac{1}{\sqrt{2}}\Big),\quad \text{if}\ \frac{1}{2}\Big(\frac{1}{2\sqrt{2}}+\frac{4}{5}\Big)t<x<0,\\
\frac{1}{2}\Big(1-2\frac{x}{t}\Big),\quad \text{if} \ 0<x<\frac{2}{5}t,\\
\frac{1}{10}, \quad \text{if}\ x>\frac{2}{5}t. 
\end{cases}
\end{equation}
In figure \ref{fig:3}, we illustrate our results.
\begin{figure}[ht!]
\begin{center}
\begin{tabular}{c c}
\includegraphics[scale=0.35]{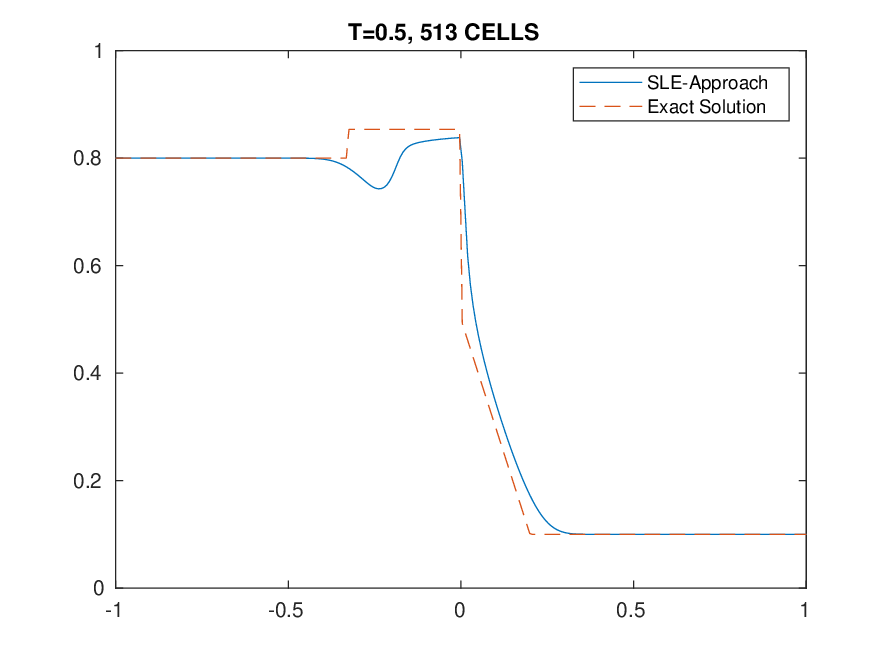} &
\includegraphics[scale=0.35]{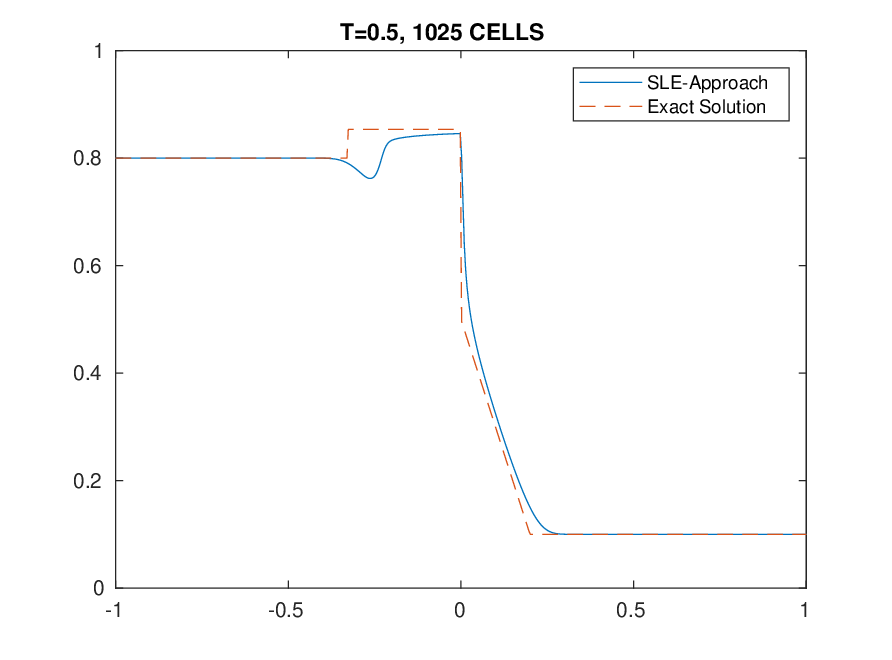} \\
($a$) & ($b$) \\
\includegraphics[scale=0.35]{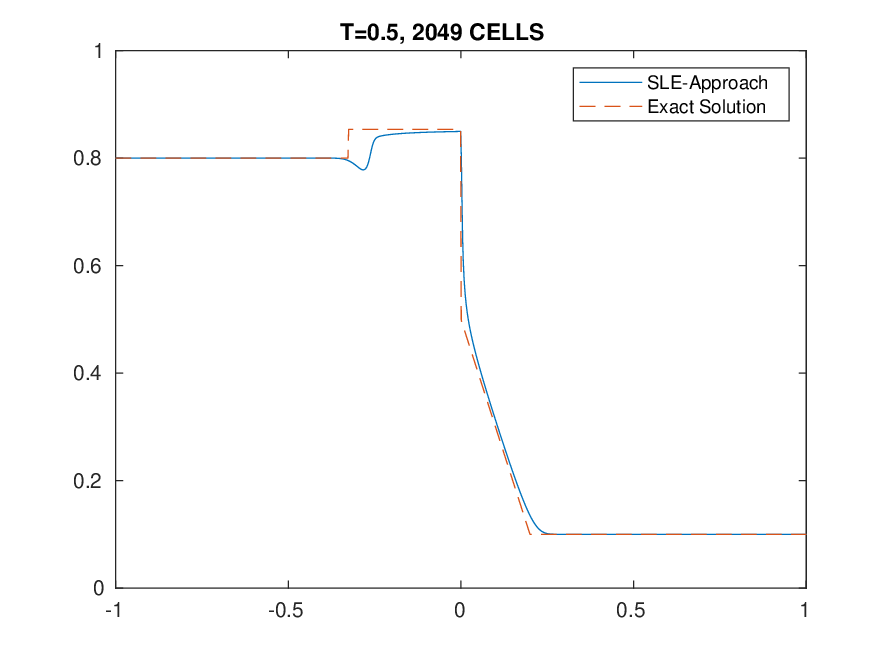} &
\includegraphics[scale=0.35]{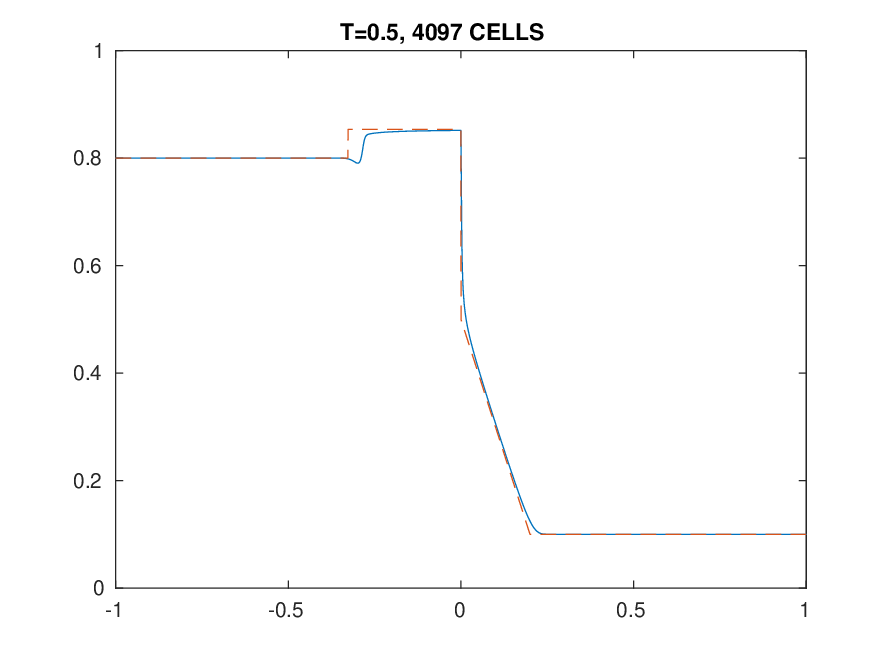} \\
($c$) & ($d$)\\
\includegraphics[scale=0.35]{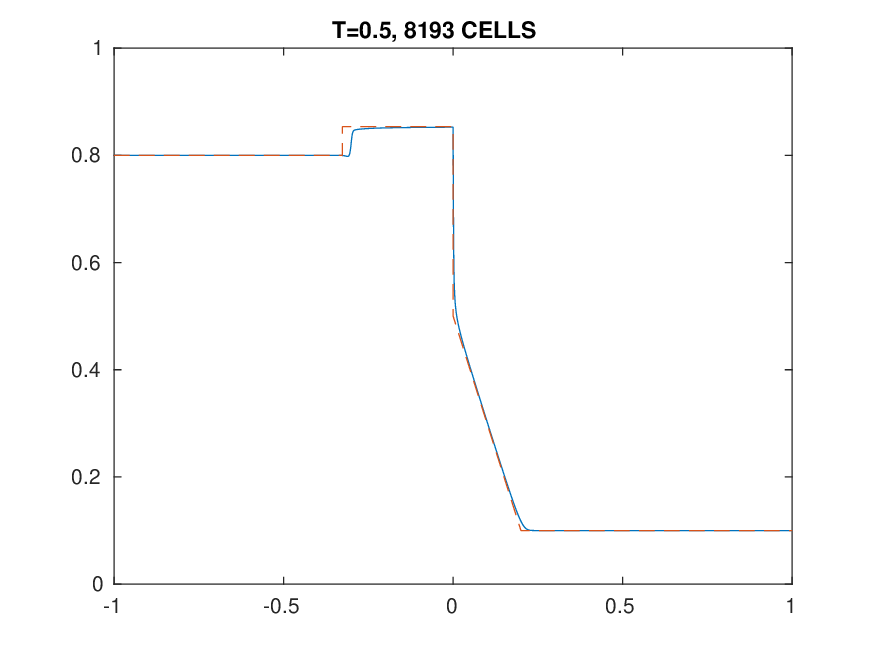} &
\includegraphics[scale=0.35]{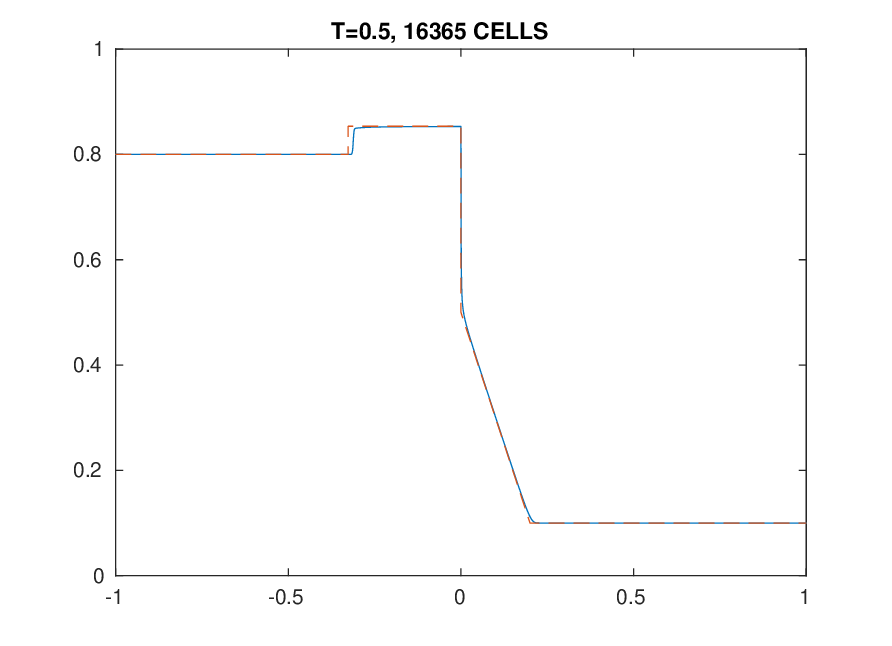} \\
($e$) & ($f$)
\end{tabular}
\end{center}
\caption{Illustration comparing the exact and approximate solutions of the Cauchy problem \eqref{1-1} with datum \eqref{Exm3-sol}. The exact solution is given by \eqref{Sol-3}. In (e) and (f), we consider a double mesh parameter to improve the structure of the solution.
}
\label{fig:3}
\end{figure}
while in Table \ref{Table:3}, we provide errors and total variation of the approximations.

\begin{table}[h!]
    \centering
\begin{tabular}{ |c|c|c|c| } 
 \hline
 cells & $L^1$-error & $L^1$-relative error &$TV$ approximate solution \\ \hline
 257 & 0.0605 & 0.0327 &0.9021 \\ 
 513 & 0.0375 &0.0202 &0.8904 \\
 1025 & 0.0225 &0.0121 &0.8672 \\ 
 2049 & 0.0130 & 0.0070&0.8428 \\
 4097 & 0.0075 &0.0040 &0.8221 \\
 8193 & 0.0042 &0.0022 & 0.8090\\
 16365 &0.0023  &0.0012 & 0.8062 \\
 
 \hline
\end{tabular}

    \caption{$L^1$-error, $L^1$-relative error and total variation of approximate  solutions for Example \ref{Example-5.8}. The computational time is 30 seconds.}
    \label{Table:3}
\end{table}

\end{exmp}
\newpage

\section{Concluding Remarks and Future Directions}\label{ConRem}
In this work, we have adapted the semi-discrete Lagrangian-Eulerian method proposed in \cite{Eduardo1} to approximate solutions for a class of scalar conservation laws with spatially discontinuous flux. Even in the case of initial data with bounded total variation, solutions to these equations may have BV norms that grow over time, and therefore the TVD property of classical schemes cannot be used. We overcome this difficulty by estimating how much the total variation of the approximate solution can grow over a time interval $[0,T]$ and applying classic compactness arguments.

Finally, we present a section of numerical experiments to illustrate the solutions obtained by our scheme, considering recent examples of scalar conservation laws with discontinuous flux in the literature. As future research directions, we plan to extend semi-discrete Lagrangian-Eulerian schemes to more general scalar conservation laws presenting multiple spatial discontinuities in the flux and investigate the compactness properties of these. In the context of traffic modeling, a natural extension will also include conservation laws on networks, where the spatial discontinuity of the flux represents an intersection of incoming and outgoing roads.

\section*{Acknowledgements}

Eduardo Abreu acknowledges the support received through the 
{\it Universidade Estadual de Campinas} (UNICAMP), 
Department of Applied Mathematics at IMECC. Richard De la cruz and 
Juan Juajibioy Otero acknowledges the support received through the 
Universidad Pedagógica y Tecnológica de Colombia (UPTC), School 
of Mathematics and Statistics. Wanderson Lambert acknowledges the 
support received through the Federal University of Alfenas UNIFAL, 
Poços de Caldas, MG, Brazil. In addition, Juan Juajibioy wants to 
express his gratitude toward Eduardo Abreu and Wanderson Lambert 
for their inspiring conversation and hospitality during the visit 
to the {\it Universidade Estadual de Campinas} (UNICAMP), Brazil.\\[0.1in]
{\bf Funding:} E. Abreu gratefully acknowledges the financial support received 
from the Brazilian National Council for Scientific and Technological 
Development (CNPq) (Grant No. 307641/2023-6) 
and the State of São Paulo Research Foundation (FAPESP) (Grant No. 2022/15108-0). Juan Carlos Juajibioy Otero acknowledges the State of São Paulo 
Research Foundation (FAPESP) from a Visiting Researcher Award 
(Grant No. 2022/15108-0).



\end{document}